\documentclass[aps,prx,preprint,superscriptaddress,amsmath,amssymb]{revtex4-2}

\usepackage{amsmath,amssymb,amsthm,mathtools}
\IfFileExists{lmodern.sty}{\usepackage{lmodern}\usepackage{microtype}}{}
\usepackage{bm}
\usepackage{graphicx,placeins}
\usepackage{booktabs,tabularx,array,enumitem,xcolor,tikz,float,needspace}
\usepackage[colorlinks=true,linkcolor=blue,citecolor=blue]{hyperref}
\graphicspath{{Figures/}}

\makeatletter
\def\switch@array{\let\@array@sw\@array@sw@LaTeX}
\makeatother

\newtheorem{theorem}{Theorem}[section]
\newtheorem{proposition}[theorem]{Proposition}
\newtheorem{lemma}[theorem]{Lemma}

\theoremstyle{definition}

\theoremstyle{remark}
\newtheorem{remark}[theorem]{Remark}

\newcommand{\M}{M}                 
\newcommand{\Mo}{M^\circ}          
\newcommand{\eps}{\varepsilon}
\newcommand{\R}{\mathbb{R}}
\newcommand{\Z}{\mathbb{Z}}
\newcommand{\E}{\mathbb{E}}
\newcommand{\one}{\mathbf 1}
\newcommand{\Exp}{\mathrm{Exp}}
\newcommand{\Ric}{\mathrm{Ric}}
\newcommand{\II}{\mathrm{II}}
\newcommand{\dv}{\operatorname{div}}
\newcommand{\avg}[1]{\langle #1\rangle}
\newcommand{\cH}{\mathcal H}
\newcommand{\cE}{\mathcal E}
\newcommand{\cS}{\mathcal S}       
\newcommand{\cZ}{\mathcal Z}       
\newcommand{\teq}{\tau_{\mathrm{eq}}} 
\DeclareMathOperator{\Var}{Var}
\DeclareMathOperator{\dist}{dist}
\DeclareMathOperator{\diam}{diam}

\begin{document}

\title{An entropic barrier to equilibration in the discrete nonlinear Schr\"odinger equation: \texorpdfstring{\\[2mm]\large}{} Rigorous bounds for the microcanonical Brownian motion}

\author{Claudio Arezzo}
\affiliation{International Centre for Theoretical Physics (ICTP), I-34151 Trieste, Italy}
\affiliation{Dipartimento di Scienze Matematiche, Fisiche e Informatiche, Università di Parma, Parco Area delle Scienze 53/A, I-43124 Parma, Italy.}

\author{Antonello Scardicchio}
\affiliation{International Centre for Theoretical Physics (ICTP), I-34151 Trieste, Italy}
\affiliation{INFN, Sezione di Trieste, I-34128 Trieste, Italy}

\date{\today}

\begin{abstract}
\noindent
The relaxation of the discrete nonlinear Schr\"odinger equation is known to slow down dramatically when the energy density $\eps$ exceeds a critical value $\eps_c$, where the system enters a {\it condensed phase}. Replacing the Hamiltonian dynamics by Brownian motion on the microcanonical surface of $N$ sites at fixed mass and energy, we prove rigorous upper and lower bounds on the equilibration time $\tau_{\rm eq}$, which is governed by the first non-zero eigenvalue $\lambda_1$ of the Laplacian on this surface. In the condensed phase the equilibration time is \emph{stretched exponential}, not exponential, in the system size. The upper bounds on $\lambda_1$ come from trial functions: of one occupation for $\eps<\eps_c$, where the problem reduces to a Sturm--Liouville equation, and of two occupations for $\eps>\eps_c$, where they cross the bottleneck through which the condensate moves between sites. The lower bounds on $\lambda_1$ come from a Bochner--Reilly argument, with a vector field that compensates the concavity of the boundary. As a result, $\tau_{\rm eq}=O(1)$ for $\eps<\eps_c$, while for $\eps>\eps_c$, asymptotically in $N$, $(\sqrt{2}-1)\sqrt{\eps-\eps_c}\leq\ln\tau_{\rm eq}/\sqrt N\leq\sqrt{\eps-\eps_c}$. We conjecture that the left-hand bound is sharp and the right-hand one is not.
\end{abstract}

\maketitle

\section{Introduction}
The success of thermodynamics in describing the macroscopic world is based on the fact that most laboratory systems achieve equilibrium states in a reasonably short time. When the microscopic laws of motion are known and trusted, time averages of observables in or about equilibrium states are described by the appropriate distributions dictated by the statistical mechanics of the corresponding microscopic problem. Such substitution is based on the {\it ergodic hypothesis}, which allows one to replace time averages of observables (to wash out microscopic fluctuations) to ensemble averages. This intuition, due to Boltzmann \cite{Boltzmann1964} for classical dynamics, and ported by Von Neumann for quantum dynamics \cite{goldstein2010normal}, served so well to understand the physical world, that it became every physicist's intuition when approaching the study of a non-exactly solvable (non-integrable) model made of many particles. It is then surprising when the ergodic hypothesis fails and, as everything which is counterintuitive, it is worth a better look because it could imply the necessity for a {\it paradigm shift} in the theory \cite{kuhn1970structure}. 

The failure modes of the ergodic hypothesis almost always involve a timescale which diverge very fast with the system size. Therefore, one can build a {\it macroscopic} timescale from a {\it microscopic} timescale and the inversion of the infinite-time and the infinite-volume limits (also known as thermodynamic limit) implicitly assumed in the ergodic hypothesis, is impossible. 

Sometimes the divergence is understood as spontaneous symmetry breaking: a ferromagnet magnetises in one direction although its Hamiltonian is rotation invariant. Sometimes the broken symmetry is not apparent in the microscopic Hamiltonian, as in mean-field spin glasses \cite{edwards1975theory,mezard1987spin} or in Anderson \cite{anderson1958absence,Evers2008} and many-body localization \cite{basko2006metal,ros2015integrals,sierant2025many}, where it is revealed only after a rewriting of the quenched average of the partition function \cite{parisi1980order,efetov1999supersymmetry}.

Ergodicity breaking therefore, for the purpose of our discussion, arises from the birth of a timescale for equilibration, $\teq$, which diverges fast, say faster than any polynomial, in the thermodynamic limit. Typically, a polynomial divergence is confined to scale-invariant critical points (critical slowing down) or is due to hydrodynamic modes (Thouless time) . Ergodicity breaking timescales are, instead, fatser than those. However, they are difficult to control quantitatively and/or analytically. For example, a divergence of the form $\teq\sim e^{cN}$ can be easily mistaken, in numerical experiments, for a large polynomial $\teq\sim N^a$, with $a\sim 10$ or for a stretched exponential $\teq\sim e^{N^\delta}$ (with $\delta<1$), although  these laws might have different physical mechanisms at their origins. A time $N^a$ with $a\gg 1$ is usually due, for example, to subdiffusion \cite{vznidarivc2016diffusive} and we do not call a subdiffusive system non-ergodic.

It is then of utmost importance to have firm, analytical results on these timescales. In this paper, we plan to do exactly that, leveraging on a series of results on the breakdown of ergodicity occurring in the discrete, non-linear Schr\"odinger equation (DNLSE) \cite{gradenigo2021localization,arezzo2022localization,iubini2025effective} and on a series of results of the eigenvalues of the Laplacian operator on large-dimensional manifolds.

The paper is organized as follows: We first present the model, the approximations and the theorems characterizing the geometry of the microcanonical surface. These are known from previous works \cite{nam2020large,gradenigo2021localization,arezzo2022localization}, and are collected in Section~\ref{sec:geometry}. The transition at $\eps=\eps_c$ is a {\it dynamical} transition \cite{gradenigo2021localization,mithun2018weakly}, where the manifold gets separated into $N$ distinct symmetry sectors, and the junctions between them have small volume and sectional area. Traversing from one sector, where a certain charge $q_i$ is large, to another, takes time which one can link to the inverse of the smallest, non-zero eigenvalue of the Laplacian on the manifold, call it $\lambda_1$. $\lambda_0=0$ is the ``trivial" equilibrium state eigenvalue.

This leads to the second part of the paper, where we estimate $\lambda_1$ from above and below. To do this we need two different techniques. Upper bounds are simpler, since they are {\it variational} in nature. We can find a set of functions on the manifold on which to calculate $\lambda_1$ and this gives a proper upper bound. The optimization in a class of problems is connected to the solution of a problem of supersymmetric quantum mechanics {\it a la} Witten \cite{witten1982supersymmetry}. We do this in Section~\ref{sec:upper}. 

Lower bounds are more complicated, since one has to find an inequality {\it valid for all functions on the manifold}. This is obtained by a completely different method which takes the name of Bochner \cite{Bochner} and Reilly \cite{Reilly}, and is explained in detail in Sections~\ref{sec:lower}--\ref{sec:corners}.

\section{The model and the results}\label{sec:model}

The DNLSE is a classical problem, where the Hamiltonian of a complex field is taken as 
\begin{equation}
    H=-\varkappa\sum_{\langle i,j\rangle\in G}(\psi_i^*\psi_j+\psi_j^*\psi_i)+\sum_{i=1}^N|\psi_i|^4,
\end{equation}
where $\langle i,j\rangle\in G$ is a graph over $N$ sites, and the canonical Poisson relations $\{\psi^*_i,\psi_j\}=i\delta_{ij}$ dictate the dynamics. Equivalently, separating $\psi_j=(x_j+ip_j)/\sqrt2$ one could write the Hamiltonian in position-momentum variables $(x,p)$, while writing in the ``charge-phase/action-angle" variables $\psi_j=\sqrt{q_j}\,e^{i\phi_j}$ we connect to the problem of Josephson junction arrays. 

The dynamics conserves both the energy and the total charge $Q=\sum_i|\psi_i|^2=\sum_i q_i$. Localization in the DNLSE consists in the following fact: For particular initial conditions, when the energy density $E/N=\eps$ is higher than a critical threshold $\eps_c$, the energy tends to concentrate in a single oscillator for a very long time \cite{eckmann2018breathers,cuneo2017energy}. Numerical calculations often fail to see the decay of such soliton \cite{mithun2018weakly,iubini2019PRL}, the question arises then, what dictates the timescale of such decay. 

If the kinetic energy term gets substituted by a random walk on the resulting microcanonical surface $M_{N,\eps}$ of constant $Q=\sum_iq_i$ and $E=\eps N=\sum_i|\psi_i|^4=\sum_iq_i^2$ and $q_i\geq 0$ (with reflecting boundary conditions and natural measure $\mu$) most of the statistical properties of the original NLSE are conserved \cite{iubini2013discrete,iubini2014coarsening}. In particular, the Brownian motion on the microcanonical surface shows the phenomenon of localization present in the original model \cite{arezzo2022localization} in all its qualitative and quantitative aspects, as explained in statistical mechanics terms \cite{gradenigo2021localization,iubini2025effective,chatterjee2014invariant}. Notice that, in this Brownian motion dynamics, the graph $G$ disappears from the problem. 

Two comments on the choice of the Brownian dynamics, which we will expand a bit in the Conclusion. The metric on $M_{N,\eps}$ is a modelling choice: other reversible dynamics with the same equilibrium measure share all the static properties discussed below, but not necessarily the subleading terms of the time scales (see the Conclusions). All the results of this paper concern this Brownian motion. Their transfer to the Hamiltonian dynamics of the DNLSE rests on this substitution, which is yet to be proven. 

In this paper we prove that the Brownian motion equilibration timescales behave in two different ways, depending we are in the normal or condensed phase, namely $\eps<2$ or $\eps>2$. First of all, however, let us define the equilibration time.

Let $P_t=e^{t\Delta_{\M}/2}$ be the semigroup of the Brownian motion on $\M_{N,\eps}$ with
reflecting conditions. $-\tfrac12\Delta_{\M}$ is self-adjoint and
non-negative on $L^2(\mu)$, with spectral gap $\lambda_1/2$, where $\lambda_1$ is the first non-zero Neumann eigenvalue of $-\Delta_\M$, eq.~\eqref{eq:lambda1} below.

The {\it relaxation} (or equilibration) {\it time} of the Brownian motion on $\M_{N,\eps}$ is
\begin{equation}
\teq=\frac2{\lambda_1}.
\label{eq:taurel}
\end{equation}
The ground state is the constant function $\rho_{\rm eq}\equiv1$ on $M$, with $\lambda_0=0$.

The relaxation time is characterised by the following, which follows from the
spectral theorem by expanding in eigenfunctions of $-\Delta_{\M}$. For every
$f\in L^2(\mu)$ and every $t\ge0$,
\begin{equation}
\bigl\|P_tf-\mu(f)\bigr\|_{L^2(\mu)}\ \le\ e^{-t/\teq}\,\bigl\|f-\mu(f)\bigr\|_{L^2(\mu)} ,
\label{eq:L2contraction}
\end{equation}
with equality when $f$ is an eigenfunction associated with $\lambda_1$. Equivalently, if the initial law has density $\rho_0$ with respect to $\mu$, the density $\rho_t$ at time $t$
satisfies
\begin{equation}
\|\rho_t-\rho_{\rm eq}\|_{L^2(\mu)}\ \le\ e^{-t/\teq}\,\|\rho_0-\rho_{\rm eq}\|_{L^2(\mu)} .
\label{eq:L2density}
\end{equation}
In particular, the equilibrium autocorrelation of any observable,
$C_f(t)=\mathbb E_\mu\bigl[(f(q_0)-\mu(f))(f(q_t)-\mu(f))\bigr]$, obeys
$|C_f(t)|\le e^{-t/\teq}\Var_\mu(f)$. The decay rate of $C_f$ equals $\lambda_1/2$
only if $f$ has a non-zero component along the eigenspace of $\lambda_1$, if not $C_f$ decays faster. The equilibration time $\teq$ is the {\it slowest} time observed in the dynamics. Some correlation functions or operator expectations will decay faster, if they are not coupled to the slowest sector, and this is observed in numerics. However, on timescales $t\gg \teq$ any signal of the initial condition will be reduced by a factor $e^{-t/\teq}$.

In this paper we prove that, for the Brownian motion defined before, there exists an $\eps_c$ ($\eps_c=2$ in our units) such that for $\eps<\eps_c$ the equilibration time is $\teq=O(1)$ and for $\eps>\eps_c$, $\teq=e^{\Theta(\sqrt N)}$. We prove an upper and a lower bound on $\lambda_1$ in each phase: the upper bounds in Section~\ref{sec:upper}, the lower bounds in Sections~\ref{sec:lower}--\ref{sec:corners}.

Without discussing the strategy, we present here the results. For $\eps<2$, and $N$ sufficiently large (i.e.\ $N>2000/(2-\eps)^2$), one has
\begin{equation}
    \frac{(2-\eps)^2}{100(\eps-1)}<\lambda_1\lesssim\lambda_{\rm MF},
\end{equation}
where $\lambda_{\rm MF}$ is the lowest non-zero eigenvalue of a second-order differential equation in one variable which can be easily tabulated to arbitrary accuracy. In the limit $\eps\to2$, $\lambda_{\rm MF}(2)=1/4$ while for $\eps\to 1$, $\lambda_{\rm MF}\simeq1/(\eps-1)$. The upper bound holds asymptotically, $\limsup_{N\to\infty}\lambda_1\le\lambda_{\rm MF}$ (Section~\ref{sec:fluid}), hence ``$\lesssim$''. The numerical experiments find a decay time which is much closer to the upper bound $2/\lambda_{\rm MF}$ than to the lower bound, which suggests that the lower bound is not strict. We could not, however, find a better lower bound.

For $\eps>2$ we have (assuming the limit exists)
\begin{equation}
(\sqrt2-1)\sqrt{\eps-2}\le\lim_{N\to\infty}\frac{\log\teq}{\sqrt N}\le \sqrt{\eps-2}. 
\end{equation}

(Theorems~\ref{thm:upper} and~\ref{thm:lower}). The stretched-exponential time scale is a consequence of the second conservation law. In condensed zero-range processes, which conserve only the mass, the condensate moves between sites on a polynomial time scale \cite{BeltranLandim}. Here it cannot move without passing through configurations in which two sites share the excess energy, whose area is smaller by a factor $e^{-(\sqrt2-1)\sqrt{N(\eps-2)}}$ (Section~\ref{sec:thermo}).

In the next Section, we introduce some results about the geometry of the microcanonical surface that will be used in Section \ref{sec:upper}. More specialized results to be used in Section \ref{sec:lower} will be introduced at the beginning of that Section.

\section{Geometry of the microcanonical surface}\label{sec:geometry}

\begin{figure}
    \centering
    \includegraphics[width=0.95\columnwidth]{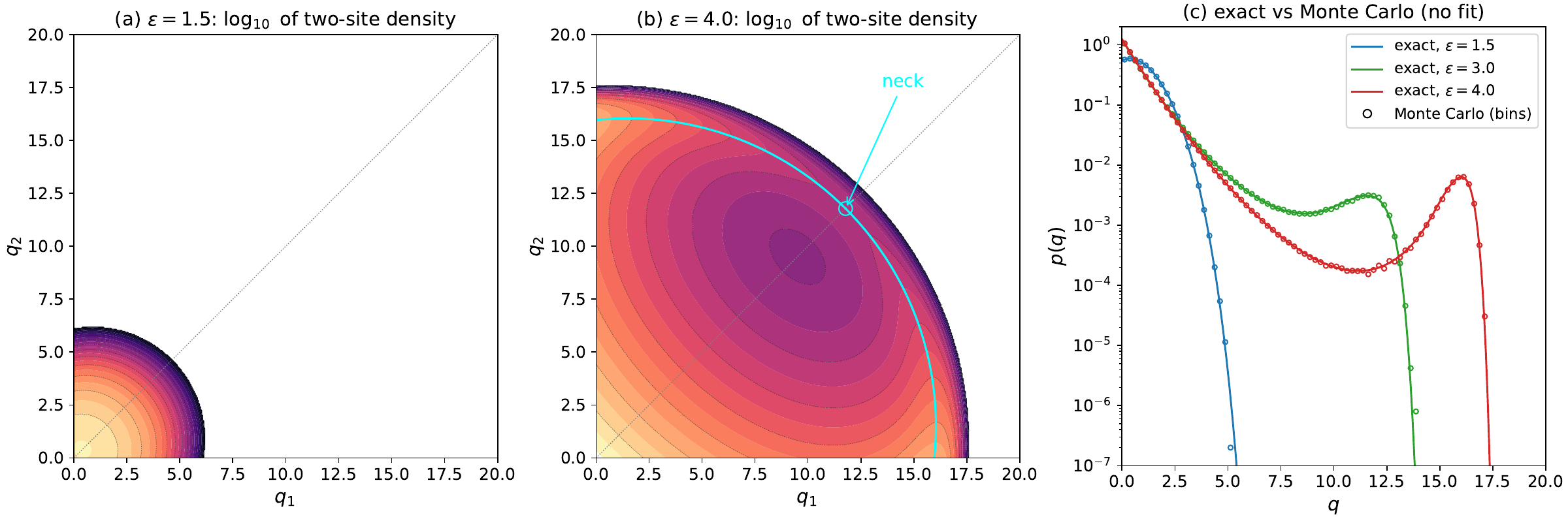}
    \caption{Several figures to understand the geometry of the high-dimensional space obtained by a numerically exact evaluation of the densities (exponential tilting plus the fluid--condensate decomposition \eqref{eq:bigjump}), at $N=100$. {\it a)} A two-dimensional projection (higher volumes are brighter colors) for $\varepsilon=1.5<2$; {\it b)} The same with $\eps=4>2$ where one can identify the two chambers $q_1\simeq a\sqrt{N}$ and $q_2\simeq a\sqrt{N}$ and, around the origin the projection of the remaining $N-2$ chambers; {\it c)} Projection on a single coordinate which shows the appearance of the chambers and the tunnel and the agreement with Montecarlo results.}
    \label{fig:sections}
\end{figure}

\begin{figure}
    \centering
    \includegraphics[width=0.95\columnwidth]{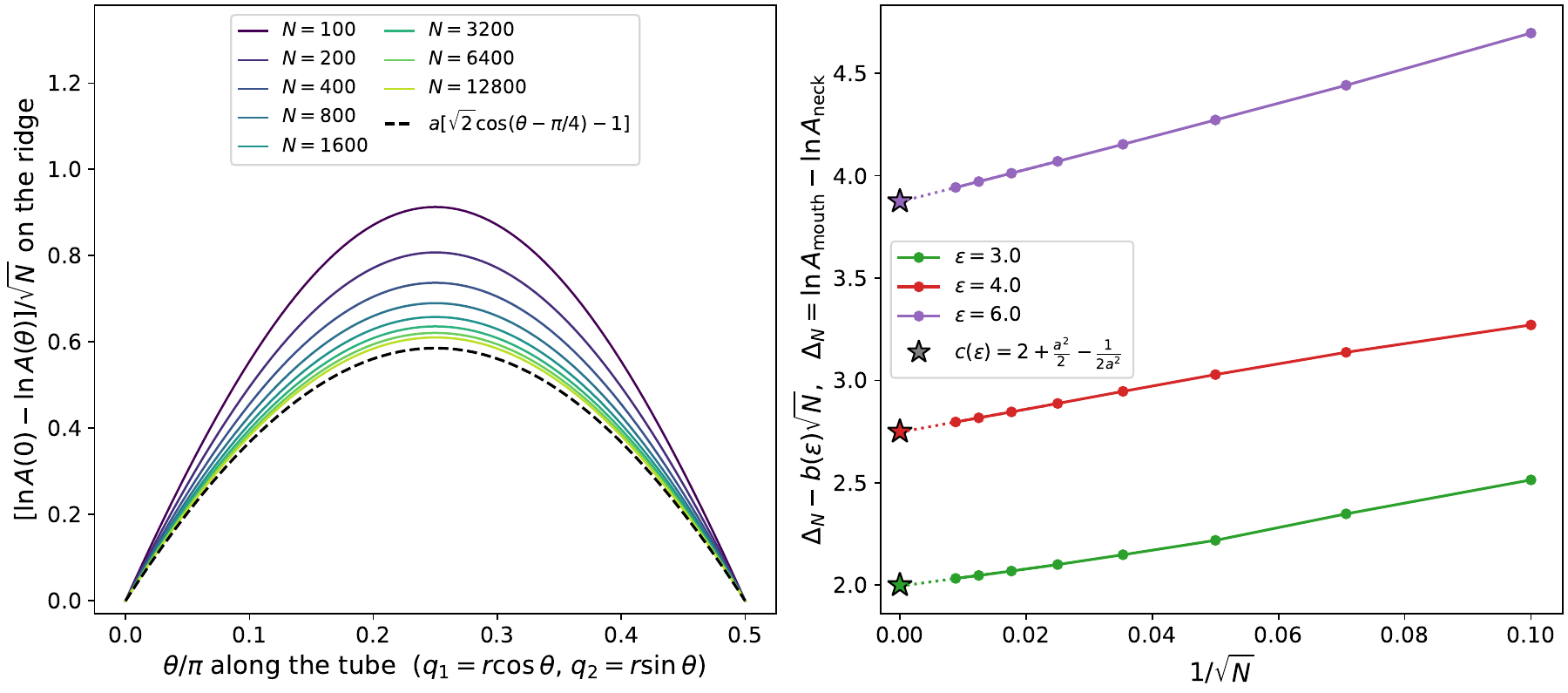}
    \caption{The numerical computed areas at $N=100,...,12800$ to show convergence to the asymptotic formula. Leading order in {\it left} and $O(1)$ corrections in {\it right} panel.}
    \label{fig:barrier}
\end{figure}

\theoremstyle{remark}
\newtheorem{elementary}[theorem]{Lemma}

\subsection{The surface and its regions}\label{sec:regions}

Throughout, $N\ge5$ and $1<\eps<N$, and we use units in which the total charge is $Q=N$. The microcanonical surface is the intersection of the simplex $\{q_i\ge0,\ \sum_iq_i=N\}$ with the sphere $\sum_iq_i^2=E=N\eps$:
\[
M=M_{N,\eps}=\Big\{q\in[0,\infty)^N:\ \sum_i q_i=N,\ \sum_i q_i^2=N\eps\Big\},\qquad R^2:=N(\eps-1).
\]
Since $|q-\one|^2=R^2$, $M$ is the part of a round sphere of radius $R$ (centre $\one$, inside the hyperplane $\sum_iq_i=N$) cut out by the $N$ half-spaces $q_j\ge0$: a spherical polytope of dimension $d=N-2\ge3$, whose boundary $\partial M=\{q\in M:\ \exists j,\ q_j=0\}$ consists of the faces $\{q_j=0\}$. We write $\sigma=\cH^{N-2}$ for its Riemannian measure, $\mu=\sigma/\sigma(M)$ for the normalised measure, and $\Mo$ for its interior. The largest value of a coordinate on $M$ is $t^*:=1+\sqrt{(N-1)(\eps-1)}$, attained only at the pole $P_i$, where $q_i=t^*$ and all the other coordinates are equal. We write $\avg{\psi}:=\frac1N\sum_j\psi(q_j)$ for a function $\psi$ of one variable, and $\II(v,w)=\langle\nabla_v\nu,w\rangle$ for the second fundamental form of a boundary face with outward normal $\nu$, so that convex faces have $\II\ge0$. For $\eps>2$ we put $a:=\sqrt{\eps-2}$.

The following are elementary standard computations which will be used throughout the paper:

\begin{elementary}\label{lem:calc}
Let $\gamma(x):=1-\frac1N-\frac{(x-1)^2}{R^2}$. On $M$, for all $i,j$:
\begin{enumerate}[label=\textup{(\arabic*)}]
\item $\langle\nabla q_i,\nabla q_j\rangle=\delta_{ij}-\frac1N-\frac{(q_i-1)(q_j-1)}{R^2}$, so that $|\nabla q_i|^2=\gamma(q_i)$; moreover $\sum_j\nabla q_j=\sum_jq_j\nabla q_j=0$.
\item $\nabla^2q_i=-\frac{q_i-1}{R^2}\,g$, $\Delta q_i=-\frac{N-2}{R^2}\,(q_i-1)$ and $\Ric=\frac{N-3}{R^2}\,g$.
\item If $\gamma(K)>0$, the face $\{q_j=K\}$ of $\{q_j\ge K\}$ has outward normal $\nu=-\nabla q_j/|\nabla q_j|$ and $\II=\frac{K-1}{R^2|\nabla q_j|}\,g$. Since $q_j=1$ at the centre $\one$, $\{q_j\ge K\}$ is larger than a hemisphere, with concave boundary, for $K<1$, and is a cap with convex boundary for $K>1$. In particular the faces $\{q_j=0\}$ of $M$ are concave.
\item For $F=\avg{\psi}$ and all $\alpha,\beta\in\R$, $|\nabla F|^2\le\frac1N\avg{(\psi'-\alpha-\beta x)^2}$; on $\{q_i=0\}$, $\langle\nabla F,\nabla q_i\rangle=\frac1N\big(\psi'(0)+\avg{(x-\eps)\psi'}/(\eps-1)\big)$.
\item $(\eps-1)\Delta F=(\eps-1)\avg{\psi''}-\avg{(x-1)\psi'}+\mathcal R_\psi$, with $|\mathcal R_\psi|\le4(\eps-1)\sup|\psi''|/N$.
\end{enumerate}
\end{elementary}

\emph{Slices and meridians.} 

Many geometric properties of $M$ will be recovered by slicing it at a given coordinate, and reducing it to a lower dimensional microcanonical surface of lower dimension and different energy. To make this idea quantitative we fix $i$, put $n=N-1$ and
\begin{gather*}
c(t)=\frac{N-t}{N-1},\qquad \rho(t)^2=R^2-\frac{N(t-1)^2}{N-1},\\
e(t)=\frac{(N-1)(N\eps-t^2)}{(N-t)^2}=1+\frac{\rho^2}{nc^2},\qquad \theta(t)=2-e(t).
\end{gather*}
For $q\in M$ with $q_i=t$, the other $n$ coordinates have sum $N-t$ and sum of squares $N\eps-t^2$, so that $y_j=q_j/c(t)$, $j\ne i$, satisfy $\frac1n\sum_{j\ne i}y_j=1$ and $\frac1n\sum_{j\ne i}y_j^2=e(t)$. The slice $\{q_i=t\}$ is thus a rescaled microcanonical surface of energy $e(t)$, and $\theta(t)$ measures how far it lies below the critical energy $2$  (to leading order $\theta(t)\approx\big(t^2-(\eps-2)N\big)/N$).

\begin{elementary}\label{lem:slices}\leavevmode
\begin{enumerate}[label=\textup{(\arabic*)}]
\item The slice $\Sigma_i(t)=M\cap\{q_i=t\}$ is $c(t)\cdot M_{n,e(t)}$, a homothetic copy of radius $\rho(t)$.
\item In the meridian coordinates $(q_j)_{j\ne i}=c(t)\one+\rho(t)\,\omega$, with $\omega$ in the unit sphere of $\{v\in\R^{n}:\sum_jv_j=0\}$, $g=dt^2/\gamma(t)+\rho^2g_S$ and $d\sigma=\rho^{N-3}\gamma^{-1/2}\,dt\,d\omega$.
\item For $t>\eps$, $e'(t)<0$ and $(c/\rho)'>0$: the admissible directions $\omega$ increase with $t$, and every meridian through $\{q_i>\eps\}$ stays in $\Mo$ and reaches the pole $P_i$.
\end{enumerate}
\end{elementary}

\begin{remark}\label{rem:margin}
Let $\eps\ge2$ and let $T_0$ be as in \eqref{eq:threshold} below. For $N$ large and $T_0\le t<t^*$, the slice $\{q_i=t\}$ lies below the critical energy $2$ with the margin
\begin{equation}\label{eq:slice-margin}
\begin{gathered}
1<e(t)<2,\qquad n\theta(t)^2\ge\max(2000,400t),\qquad \theta(t)^2\ge\kappa_N,\\
\kappa_N=\begin{cases}
576\sqrt{\eps-2}\,N^{-1/2},&\eps>2,\\
36864N^{-2/3},&\eps=2,
\end{cases}
\end{gathered}
\end{equation}
as follows from $\theta(t)(N-t)^2=(N+1)t^2-4Nt-(\eps-2)N^2+N\eps$.
\end{remark}

\begin{figure}[tbp]
\centering
\includegraphics[width=0.985\columnwidth]{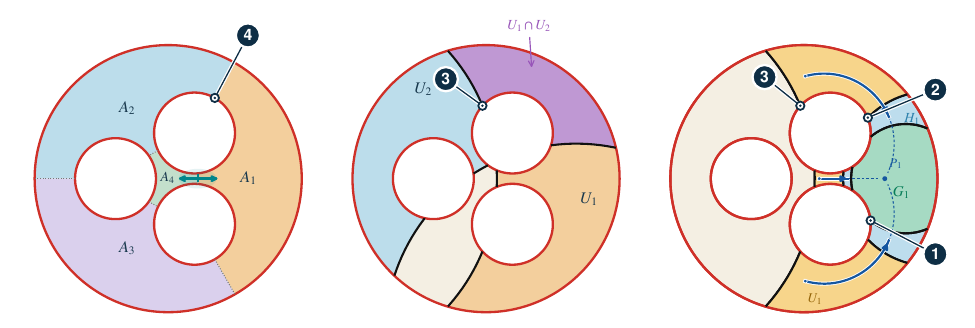}\\[2mm]
\includegraphics[width=0.985\columnwidth]{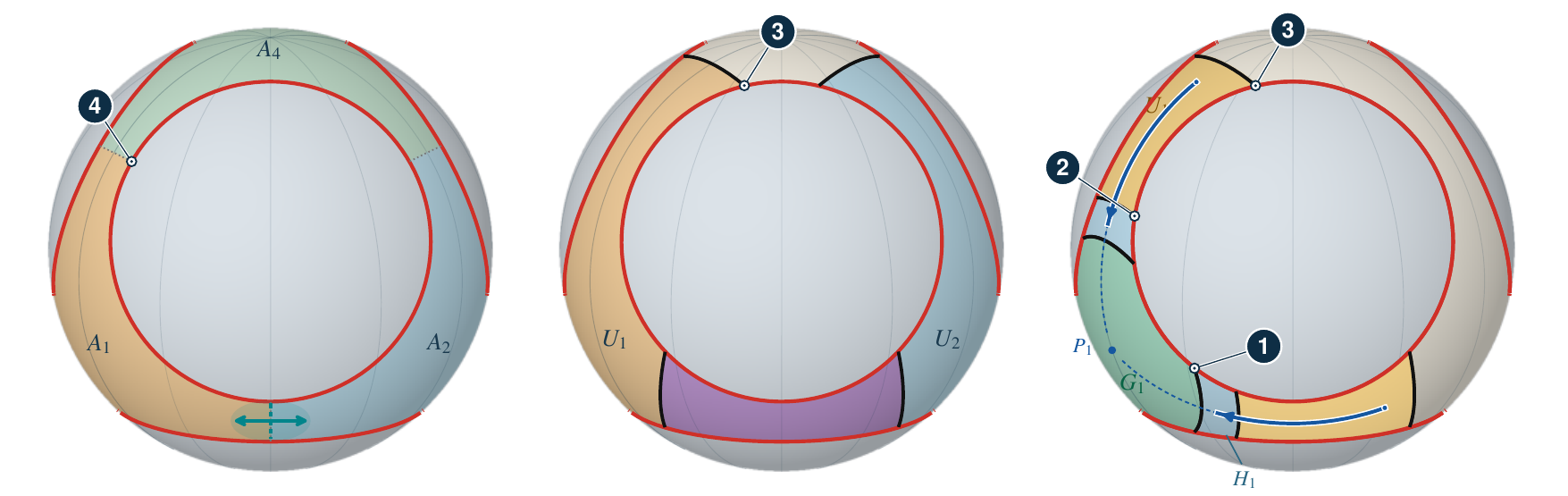}
\caption{The regions $G_i\subset H_i\subset U_i$ and the cells $A_i$ of $M$, in a planar view and on the sphere.}
\label{fig:corners-regions}
\end{figure}

\emph{The regions.} For $K\ge\eps$, the cap $\{q\in M:\ q_i>K\}$ is connected, since its meridians reach $P_i$. Its boundary consists of the face $\{q_i=K\}$, which is convex because $K>1$, and of pieces of the concave faces $\{q_j=0\}$, $j\ne i$. We will use the cells $A_i$ and three such caps,
\begin{gather*}
A_i=\{q\in\Mo:\ q_i>q_j\ \forall j\ne i\},\qquad
U_i=\{q\in\Mo:\ q_i>\eps\},\\
H_i=\{q\in\Mo:\ q_i>2\eps\},\qquad
G_i=\{q\in M:\ q_i\ge T\},
\end{gather*}
where the level of the deep cap $G_i$, used only for $\eps\ge2$, is
\begin{equation}\label{eq:threshold}
T\in[T_0,T_0+1],\qquad
T_0=\begin{cases}
\sqrt{\eps-2}\,\sqrt N+16(\eps-2)^{-1/4}N^{1/4},&\eps>2,\\
16N^{1/3},&\eps=2,
\end{cases}
\end{equation}
chosen outside the finitely many levels at which the face $\{q_i=T\}$ does not meet the faces $\{q_j=0\}$ transversally (Proposition~\ref{prop:corners}(b)). For $N$ large, $2\eps<T<T+1<t^*\le N/2$ (the last inequality holds, for example, when $N\ge4\eps$), so that $\emptyset\ne G_i\subset H_i\subset U_i$. The regions are shown in Figure~\ref{fig:corners-regions}.

\begin{elementary}\label{lem:cover}
Let $1<\eps<N/2$. The largest coordinate of every $q\in\Mo$ exceeds $\eps$, because $\eps=\sum_j(q_j/N)\,q_j$ is a weighted average of the coordinates. Hence $\Mo=\bigcup_iU_i$ and $A_i\subset U_i$; the cells $A_i$ partition $M$ up to null sets, and $\mu(A_i)=1/N$. Moreover $U_i\cap U_j\ne\emptyset$, since it contains the configurations with $q_i=q_j>\eps$ and all other coordinates equal, so that $M$ is connected.
\end{elementary}
The connectedness of $M$ for $1<\eps<N/2$ is due to \cite[Sec.~6 and App.~A]{arezzo2022localization}.

\emph{Neumann eigenvalues.} A probability density $\rho_t$ on $M$ evolves by $\partial_t\rho_t=\frac12\Delta_{M}\rho_t$, with Neumann (reflecting) conditions on $\partial M$. The measure $\mu$ is reversible, and the relaxation time \eqref{eq:taurel} is $\teq=2/\lambda_1$, where
\begin{equation}
\lambda_1=\lambda_1(N,\eps)=\inf\Big\{\frac{\int_M|\nabla f|^2\,d\mu}{\Var_\mu f}:\ f\in H^1(M),\ \Var_\mu f>0\Big\}
\label{eq:lambda1}
\end{equation}
is the first nonzero Neumann eigenvalue of $-\Delta_M$. We write $\cE(f,f)=\int_M|\nabla f|^2d\mu$ for the Dirichlet form. For a subdomain $D\subset M$, $\lambda_1(D)$ is the first nonzero Neumann eigenvalue of $D$, and $P(D)=1/\lambda_1(D)$ is its Poincar\'e constant, the smallest constant with $\min_{c\in\R}\int_D|f-c|^2\,d\sigma\le P(D)\int_D|\nabla f|^2\,d\sigma$ for all $f\in H^1(D)$. We abbreviate $P_G=P(G_i)$, $P_H=P(H_i)$, $P_U=P(U_i)$.

\subsection{Areas}\label{sec:measure}

The areas of $M$ and of its coordinate-defined subsets are expressed through the microcanonical partition function
\begin{equation}
\Omega_n(S,E)=\int_{[0,\infty)^n}\!\!d^nq\;
\delta\Bigl(S-\sum_i q_i\Bigr)\,\delta\Bigl(E-\sum_i q_i^2\Bigr).
\end{equation}
The constraint gradients $\nabla\sum_iq_i=\one$ and $\nabla\sum_iq_i^2=2q$ have the Gram matrix
\begin{equation}
G=\begin{pmatrix} N & 2N\\[2pt] 2N & 4N\eps\end{pmatrix},\qquad \sqrt{\det G}=2N\sqrt{\eps-1},
\label{eq:gram}
\end{equation}
which is constant on $M$. Moreover, by the coarea formula:
\begin{equation}
d\mu(q)=\frac{1}{\Omega_N(N,N\eps)}\;
\delta\Bigl(N-\sum_i q_i\Bigr)\,\delta\Bigl(N\eps-\sum_i q_i^2\Bigr)
\prod_{i=1}^{N}\mathbf 1_{\{q_i\ge0\}}\,dq_i ,
\label{eq:muconstraint}
\end{equation}
and, since $e^{-\sum_iq_i}=e^{-N}$ on $M$, the conditioned product of exponentials
\begin{equation}
\mu=\mathrm{Law}\Bigl(X\ \Big|\ \sum_i X_i=N,\ \sum_i X_i^2=N\eps\Bigr),
\qquad X_i\ \text{i.i.d.}\ \Exp(1).
\label{eq:muprob}
\end{equation}
For any number $n$ of coordinates, let $X_1,\dots,X_n$ be i.i.d.\ $\Exp(1)$. For $L\in(0,\infty]$ let $f_n^{<L}(S,E)$ be the density of $(\sum_jX_j,\sum_jX_j^2)$ restricted to the event $\{\max_jX_j<L\}$, and $f_n:=f_n^{<\infty}$; let $\Omega^{<L}_n$ be defined as $\Omega_n$, with all $n$ coordinates restricted to $[0,L)$.

\begin{elementary}[areas as densities]\label{lem:density}\leavevmode
\begin{enumerate}[label=\textup{(\arabic*)}]
\item $f_n^{<L}(S,E)=e^{-S}\,\Omega_n^{<L}(S,E)$. More generally, if $X_1,\dots,X_n$ are i.i.d.\ with a density $p$, the density of $(\sum_jX_j,\sum_jX_j^2)$ is bounded by $\frac\pi{\sqrt3}(\sup p)^3$ for $n\ge3$; in particular $f_n^{<L}\le\pi/\sqrt3$.
\item $\sigma(M_{N,\eps})=2N\sqrt{\eps-1}\,\Omega_N(N,N\eps)=2N\sqrt{\eps-1}\,e^N\,f_N(N,N\eps)$.
\item For sets defined by conditions on $k\le N-3$ coordinates, areas are obtained by integrating out those coordinates. For instance the joint density of two coordinates is
\[
\mu\big(q_i\in du,\ q_j\in dv\big)=\frac{\Omega_{N-2}\big(N-u-v,\ N\eps-u^2-v^2\big)}{\Omega_N(N,N\eps)}\,du\,dv,
\]
so that
\[
\mu(U_i\cap U_j)=\frac1{f_N(N,N\eps)}\int_\eps^\infty\!\!\int_\eps^\infty e^{-u-v}\,f_{N-2}\big(N-u-v,\ N\eps-u^2-v^2\big)\,du\,dv .
\]
\end{enumerate}
\end{elementary}

\begin{elementary}[local limit theorem {\cite[\S19]{BR}}]\label{lem:llt}
Uniformly in $z\in\R^2$, $n\,f_n(n+\sqrt n z_1,\,2n+\sqrt nz_2)$ converges to a nondegenerate Gaussian density in $z$. In particular $f_n\ge c/n$ within $O(\sqrt n)$ of $(n,2n)$.
\end{elementary}

Lemma~\ref{lem:llt} controls the partition functions $\Omega_n(S,E)=e^Sf_n(S,E)$ only within $O(\sqrt n)$ of $(n,2n)$. Far from it, at energy densities $nE/S^2>2$, where the residual system is supercritical, we use a Chernoff bound with a cutoff: when all coordinates lie below a cutoff, the energy fugacity may be taken negative.

\begin{lemma}[tilted bound]\label{lem:tilt}
Fix $\eps>2$, $a=\sqrt{\eps-2}$ and $k\in\{1,2\}$, and let $n=N-k$. Uniformly in $\varsigma\ge0$, $r\ge0$ and $\zeta\in[1/(2\sqrt N),\sqrt\eps\,]$,
\begin{equation}
\ln\Omega^{<\zeta\sqrt N}_n\bigl(N-\varsigma\sqrt N,\ N(\eps-r^2)\bigr)\ \le\ N-\sqrt N\,\Big[\varsigma+\frac{(a^2-r^2)_+}{\max(\zeta,a/2)}\Big]+C\ln N ,
\label{eq:tiltbound}
\end{equation}
with $C$ depending only on $\eps$, and $(x)_+=\max(x,0)$.
\end{lemma}
The proof is in Appendix~\ref{app:tilt}.

For $\eps>2$ the area of $M$ is dominated by the configurations in which one coordinate, of size about $a\sqrt N$, carries the whole excess energy $a^2N$ (Section~\ref{sec:thermo}). The following two Propositions give the fundamental estimates on volumes needed throughout our paper and will be proved in Appendix~\ref{app:Ceps}:

\begin{proposition}[total area]\label{prop:total}
For fixed $\eps>2$ and $N$ large,
\begin{equation}
\exp\bigl\{N-a\sqrt N-O(\ln N)\bigr\}\ \le\ \Omega_N(N,N\eps)\ \le\ \exp\bigl\{N-a\sqrt N+O(\ln N)\bigr\},
\label{eq:normlower}
\end{equation}
that is, $e^{-a\sqrt N-O(\ln N)}\le f_N(N,N\eps)\le e^{-a\sqrt N+O(\ln N)}$.
\end{proposition}
The gluing in the lower bound needs the caps $U_i$ to overlap substantially. By Lemma~\ref{lem:density}(3) the overlap is a ratio of values of the densities $f_n$; the underlying phenomenon, condensation on the intersection of a simplex and a sphere, is described in Section~\ref{sec:thermo} and in \cite{Chatterjee,nam2020large}.

Combined with the local limit theorem, the upper bound in Proposition~\ref{prop:total} has the following consequence, which is crucial for the proof of the lower bound when $\eps\ge2$; part~(a), at $\eps=2$, needs only the local limit theorem:
\begin{proposition}\label{prop:overlap}
Let $i\ne j$.
\begin{enumerate}[label=\textup{(\alph*)}]
\item At $\eps=2$, $\mu(U_i\cap U_j)\ge\alpha_2:=\frac12(e^{-3}-e^{-4})^2$ for $N\ge N_0$.
\item For fixed $\eps>2$ there is $C_\eps$ such that $\mu(U_i\cap U_j)\ge N^{-C_\eps}$ for $N\ge N_\eps$; any $C_\eps>2$ is admissible.
\end{enumerate}
\end{proposition}

Finally, along the meridians of a cap: by Lemma~\ref{lem:slices}(2), on $\{q_i>\eps\}$ one has $d\sigma=\upsilon(t)\,dt\,d\omega$, with $\upsilon(t)=\rho(t)^{N-3}\gamma(t)^{-1/2}$ and $\omega$ ranging over the admissible directions. For $t,s\in(\eps,2\eps+1]$ and $N$ large,
\begin{equation}\label{eq:meridian-density}
\gamma(t)\ge\tfrac12,\qquad \frac{\upsilon(t)}{\upsilon(s)}\le C_\eps .
\end{equation}
Deeper in the cap the area decays like $e^{-q_i}$; in the chart of $H_i$ this is quantified by Lemma~\ref{lem:chart}(a).

\subsection{The landscape: chambers and the tube}\label{sec:thermo}

This subsection explains where the two proofs come from. For $\eps>2$ the surface consists of $N$ chambers, one for each site that can carry the condensate. They are joined pairwise by tubes whose narrowest section, the \emph{neck}, is smaller than a chamber by a factor $e^{-b(\eps)\sqrt N}$, and the spectral gap is governed by this neck.
\begin{itemize}
\item \emph{Upper bound (Section~\ref{sec:upper}).} The landscape dictates the test function: a function of two occupations, $+1$ when the condensate sits on site $1$, $-1$ when it sits on site $2$ and $0$ otherwise, which changes at the neck level $q\simeq a\sqrt{N/2}$, so that its Dirichlet energy is set by the area of the necks. Proposition~\ref{prop:valley} makes the neck estimate rigorous.
\item \emph{Lower bound (Sections~\ref{sec:lower}--\ref{sec:corners}).} The landscape explains why the Bochner argument must be localised to the deep caps $G_i$, which lie inside a single chamber, and why the lower bound loses against the neck (Remark~\ref{rem:exponent}).
\end{itemize}
The subsection itself is heuristic. Apart from the exact identities \eqref{eq:rescaled} and \eqref{eq:bigjump}, it uses the asymptotics of $\Omega_n$ from \cite{gradenigo2021localization}, and its displays hold to the indicated orders without error control. It describes the landscape behind both proofs, but none of the proofs uses it: the neck exponent $b(\eps)$ of \eqref{eq:tubearea} is proved as an upper bound in Proposition~\ref{prop:valley}, and the total area as Proposition~\ref{prop:total}.

\emph{Entropy.} The computation of $\Omega_N$ has been done in \cite{gradenigo2021localization}. By homogeneity it suffices to know $\Omega_n(n,nx)$ as a function of the energy density $x$. With fugacities $\eta$ and $\beta$ conjugate to mass and energy, and the single-site integral
\begin{equation}
\cZ(\eta,\beta)=\int_0^\infty dq\;e^{-\eta q-\beta q^2},
\label{eq:gsingle}
\end{equation}
finite for $\beta\ge0$, the saddle-point equations $\E_p[q]=1$, $\E_p[q^2]=x$ for the tilted density $p(q)=e^{-\eta q-\beta q^2}/\cZ(\eta,\beta)$ have a solution with $\beta>0$ when $1<x<2$, and
\begin{equation}
\ln\Omega_n(n,nx)=n\,\cS(x)+O(\ln n),\qquad
\cS(x)=\eta+\beta x+\ln \cZ(\eta,\beta).
\label{eq:sdef}
\end{equation}
$\cS$ is concave, with maximum $\cS(2)=1$ at the infinite-temperature point $x=2$ ($\beta=0$, $\eta=1$, $p(q)=e^{-q}$), and $\cS(x)\simeq1-(x-2)^2/8$ there. As $x\to1^+$, $\beta\to+\infty$; for $1<\eps\le N/(N-1)$ the surface is a round sphere of radius $R$ inside the simplex, and $\lambda_1=(N-2)/R^2$ exactly (compare Theorem~\ref{thm:lower}(i)).

\emph{The fluid phase $1<\eps<2$.} Locally $\mu$ is a product of the densities $p(q)$, which decay faster than $e^{-q}$: the largest occupation is $O(\ln N)$ and there are no chambers. Since $\cS$ is increasing for $x<2$, an excursion $q_i=z_i\sqrt N$ of one or two sites costs a factor $e^{-N[\cS(\eps)-\cS(\eps-z_1^2-z_2^2)]}$. There is no tunnel to cross, consistent with Theorem~\ref{thm:lower}(i). For $\eps>2$ the saddle-point equations have no solution with $\beta\ge0$: the saddle sticks at the branch point $\beta=0$, the residual system stays critical, and the excess energy $N(\eps-2)$ is carried by a single site, the \emph{condensate}.

\emph{Two-site sections for $\eps>2$.} A passage from one chamber to another involves two coordinates, the one losing the excess energy and the one gaining it. We use the two-site density of Lemma~\ref{lem:density}(3) with $q_i=z_i\sqrt N$. The residual system has $n=N-2$ sites, mass $S'=N-q_1-q_2$ and energy $E'=N\eps-q_1^2-q_2^2$; the exact homogeneity $\Omega_n(S,E)=m^{n-3}\Omega_n(S/m,E/m^2)$ with $m=S'/n$ gives
\begin{equation}
\ln\Omega_{n}(S',E')=(n-3)\ln m+\ln\Omega_n(n,nx),
\qquad
x=\frac{E'n}{S'^2},
\label{eq:rescaled}
\end{equation}
where
\begin{equation}
(n-3)\ln m=-(z_1+z_2)\sqrt N-\frac{(z_1+z_2)^2}{2}+2+O(N^{-1/2}),
\label{eq:massdepletion}
\end{equation}
\begin{equation}
x=\bigl(\eps-z_1^2-z_2^2\bigr)\Bigl[1+\frac{2(z_1+z_2)}{\sqrt N}\Bigr]+O(N^{-1}).
\label{eq:xexpansion}
\end{equation}
The two spikes remove energy $q_i^2$ but also mass $q_i$; the mass removal is the $O(N^{-1/2})$ correction in \eqref{eq:xexpansion}, which displaces the tunnel by $O(1)$. Since $\ln\Omega_n(n,nx)=n-\sqrt{n(x-2)}+O(\log n)$ for $x>2$ \cite{gradenigo2021localization},
\begin{equation}
\ln\Omega_{N-2}\bigl(N-q_1-q_2,\,N\eps-q_1^2-q_2^2\bigr)=N-\sqrt N\,\Phi(z_1,z_2)+O(\ln N),
\label{eq:Phi2}
\end{equation}
with $\Phi(z_1,z_2)=z_1+z_2+\sqrt{a^2-z_1^2-z_2^2}$ for $z_1^2+z_2^2\le a^2$. $\Phi$ is symmetric, since transposing two coordinates is an isometry of $M$ that preserves $\mu$. At $z=0$ it gives $\ln\Omega_N(N,N\eps)=N-a\sqrt N\equiv\ln V_c$, the volume of the surface and, to exponential order, of a single chamber. The transverse section at fixed $(q_1,q_2)$ is therefore
\begin{equation}
\frac{\mathcal A(z_1,z_2)}{V_c}=\exp\Bigl\{-\sqrt N\bigl[\Phi(z_1,z_2)-a\bigr]\Bigr\},
\qquad V_c=e^{\,N-a\sqrt N}.
\label{eq:areaprofile}
\end{equation}

\emph{The tube.} The square root in \eqref{eq:Phi2} vanishes on the circle
\begin{equation}
z_1^2+z_2^2=a^2 ,
\label{eq:budget}
\end{equation}
where the two peeled coordinates carry the whole excess energy and the other $N-2$ are exactly critical. Inside the disk a third coordinate hosts a spike of size $\sqrt{(a^2-|z|^2)N}$, so a path through the interior is a concatenation of two passages and is not cheaper. The direct passage between two chambers is therefore confined to \eqref{eq:budget}, where $\Phi=\|z\|_1$ at fixed $\|z\|_2=a$. With $z=a(\cos\theta,\sin\theta)$,
\[
\Phi-a=a\Bigl[\sqrt2\cos\bigl(\theta-\tfrac\pi4\bigr)-1\Bigr],
\]
which vanishes at the ends $\theta=0,\pi/2$ (the chambers) and is maximal at $\theta=\pi/4$. At finite $N$ the mass removed by the peeled coordinates pushes the tube outwards (Figure~\ref{fig:sections}, panel {\it b}). Balancing the Gaussian fluctuation of the residual system, $\ln\Omega_n(n,nx)\simeq\ln\Omega_n(n,2n)-n(x-2)^2/8$ for $|x-2|\ll n^{-1/3}$, against the linear cost in \eqref{eq:massdepletion} places the tube at
\begin{equation}
\|z\|_2=a+\Bigl(2-\frac1{a^2}\Bigr)\frac{\sqrt2\cos(\theta-\frac\pi4)}{\sqrt N}+O(N^{-1}),
\label{eq:tuberadius}
\end{equation}
each large coordinate being shifted by $2-1/a^2$. Collecting the $O(1)$ terms of \eqref{eq:massdepletion}, \eqref{eq:tuberadius} and of the Gaussian fluctuation, the section along the tube is
\begin{equation}
\ln\frac{\mathcal A(\theta)}{V_c}=-\sqrt N\,a\Bigl[\sqrt2\cos\bigl(\theta-\tfrac\pi4\bigr)-1\Bigr]
-c_*(\eps)\sin2\theta+O(N^{-1/2}),
\qquad c_*(\eps)=2+\frac{a^2}{2}-\frac1{2a^2},
\label{eq:tubeO1}
\end{equation}
valid for $a\gg N^{-1/6}$; a numerical evaluation of $\Omega_n$ up to $N=12800$ confirms \eqref{eq:tuberadius} and the value of $c_*(\eps)$. In the coordinates $\bar q=(q_1+q_2)/2$, $\tau=q_1-q_2$ one has $\sin2\theta=1-\tau^2/(a^2N)$ on the tube, and the section of the tunnel joining two chambers is
\begin{equation}
\mathcal A(\tau)=V_c\,
\exp\Bigl\{-\sqrt N\Bigl[\sqrt{2a^2-\tau^2/N}-a\Bigr]-c_*(\eps)\Bigl(1-\frac{\tau^2}{a^2N}\Bigr)\Bigr\},
\qquad |\tau|\le a\sqrt N,
\label{eq:tubearea}
\end{equation}

\emph{The neck.} The narrowest point is $\tau=0$, where $q_1=q_2=a\sqrt{N/2}+2-1/a^2$ and
\[
\mathcal{A}(0)=V_{c}\,e^{-\sqrt{N}\,b(\eps)-c_*(\eps)},
\qquad b(\eps)=(\sqrt2-1)\,a .
\]
It is exactly the midpoint of the tube: at $q_1=a\sqrt{N/2}$ the excess $Na^2/2$ left to the other sites must be carried by a single coordinate, so $q_2=q_1$. Expanding the square root,
\begin{equation}
\ln\frac{\mathcal A(\tau)}{\mathcal A(0)}=\frac{\tau^2}{2\sqrt2\,a\sqrt N}
+c_*(\eps)\frac{\tau^2}{a^2N}
+O\!\left(\frac{\tau^4}{N^{3/2}}\right),
\label{eq:neckwidth}
\end{equation}
so the section stays within a constant factor of its minimum over $|\tau|\lesssim(2\sqrt2\,a)^{1/2}N^{1/4}$. The bottleneck is a window of width $N^{1/4}$ in a channel of length $2a\sqrt N$.

Figures~\ref{fig:sections} and~\ref{fig:barrier} are obtained from a numerically exact evaluation of $\Omega_n(S,E)$, by exponential tilting and a decomposition into fluid and condensate (Appendix~\ref{app:numerics}); it agrees with an exact-constraint Montecarlo at $N=100$.

\section{\texorpdfstring{Upper bound on $\lambda_1$}{Upper bound on lambda1}}
\label{sec:upper}

One can obtain upper bounds easily by variational methods. By \eqref{eq:lambda1}, every $f$ with $\E_\mu f=0$ gives $\lambda_1\le\E_\mu|\nabla_\M f|^2/\E_\mu f^2$, and the closer $f$ is to the eigenfunction, the better the bound. The natural trial functions are different in the two phases, in particular the number of variables they must depend non-trivially upon changes from one for $\eps\le 2$ to two for $\eps\geq 2$.

\subsection{\texorpdfstring{The fluid phase: $\eps<2$}{The fluid phase: eps<2}}\label{sec:fluid}

In the fluid phase, $\eps<2$, all the variables have the same order of magnitude $q_i\sim O(1)$ and the natural trial functions are functions of a single variable only, $f(q)=\varphi(q_1)-\E_\mu\varphi(q_1)$. In the limit $N\to\infty$ the
sites become independent, and the slowest mode of a product of independent processes is a
function of one coordinate. For such functions both averages reduce to one-dimensional
integrals. Since $\nabla_{\M}f=\varphi'(q_1)\,\nabla_{\M}q_1$, and $|\nabla_{\M}q_1|^2=\gamma(q_1)=1-\tfrac{1}{N}-\tfrac{(q_1-1)^2}{N(\eps-1)}$ is a function of $q_1$ alone (Lemma~\ref{lem:calc}), one finds exactly
\begin{equation}
\lambda_1\ \le\ \frac{\int_0^\infty dq_1\;\gamma(q_1)\,\varphi'(q_1)^2\,p_N(q_1)}
{\int_0^\infty dq_1\;\bigl(\varphi(q_1)-\bar\varphi\bigr)^2\,p_N(q_1)},
\qquad
p_N(q_1)=\frac{\Omega_{N-1}(N-q_1,\,N\eps-q_1^2)}{\Omega_N(N,N\eps)} ,
\end{equation}
where $p_N$ is the marginal of one occupation and $\bar\varphi=\int_0^\infty\varphi(q)\,p_N(q)\,dq$. First of all, since $\gamma(q_1)=1-O(1/N)<1$ (and goes to 1 when $N\to\infty$, for any $\eps>1$) we can replace this upper bound with another, simpler one in which $\gamma\to 1$:
\begin{equation}
\lambda_1\ \le\ \frac{\int_0^\infty dq_1\;\,\varphi'(q_1)^2\,p_N(q_1)}
{\int_0^\infty dq_1\;\bigl(\varphi(q_1)-\bar\varphi\bigr)^2\,p_N(q_1)}.
\label{eq:onesiteRQ}
\end{equation}
The best choice of $\varphi$ solves the Euler--Lagrange equation of \eqref{eq:onesiteRQ}, a one-dimensional Sturm--Liouville problem,
\begin{equation}
-\frac{1}{p_N}\frac{d}{dq}\Bigl(p_N\,\frac{d\varphi}{dq}\Bigr)=\Lambda_N\,\varphi ,
\qquad \varphi'=0\ \text{at the ends of the support},
\label{eq:SLfinite}
\end{equation}
and its first non-zero eigenvalue $\Lambda_N(\eps)$ is an upper bound on $\lambda_1$ at every finite $N$.

\begin{figure}
    \centering
    \includegraphics[width=0.95\columnwidth]{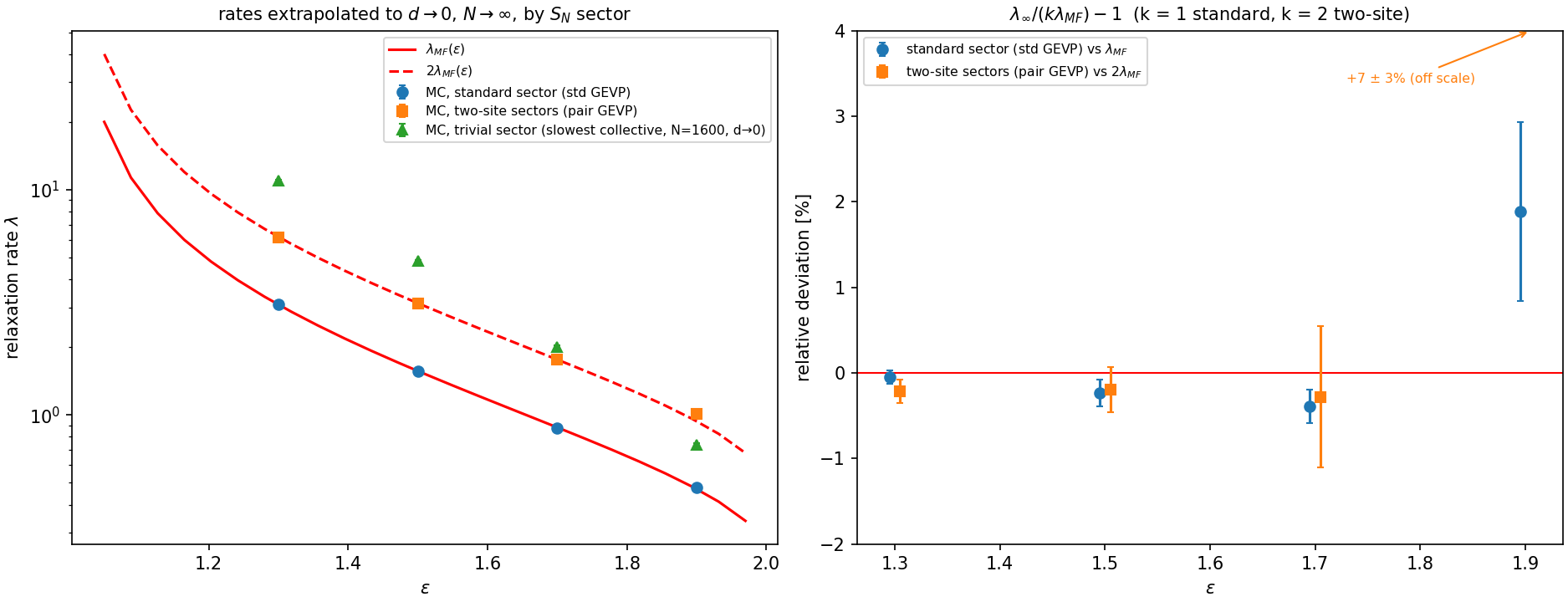}
    \caption{The analytic value of the upper bound $\lambda_{\rm MF}$ compared with the numerical values extracted from the correlation functions of the Montecarlo procedure discussed in the text. The deviations at $\eps=1.9$ are indicative of the critical region.}
    \label{fig:lmf_lMC}
\end{figure}

By the equivalence of ensembles (Section~\ref{sec:thermo}, see also \cite{nam2020large}) $p_N(q)$ tends to the tilted density $p(q)=e^{-\eta q-\beta q^2}/\cZ(\eta,\beta)$ of \eqref{eq:gsingle}, with $(\eta,\beta)$ fixed by $\E_p[q]=1$ and $\E_p[q^2]=\eps$. Equation~\eqref{eq:SLfinite} then
becomes the eigenvalue problem of a single reflected diffusion in the potential
$W(q)=\eta q+\beta q^2$,
\begin{equation}
-\varphi''+(\eta+2\beta q)\,\varphi'=\lambda\,\varphi ,\qquad q\ge0,\qquad \varphi'(0)=0 .
\label{eq:onesiteODE}
\end{equation}
This equation is solved exactly. In the variable $\xi=\sqrt{2\beta}\,q+\xi_0$, with
$\xi_0=\eta/\sqrt{2\beta}$, it becomes the Hermite equation $\varphi''-\xi\varphi'+\nu\varphi=0$ with
$\nu=\lambda/2\beta$. Its solution with polynomial growth is $\varphi=He_\nu(\xi)$, and the Neumann
condition at $\xi=\xi_0$ selects $D_{\nu-1}(\xi_0)=0$, where $D_\nu$ is the parabolic cylinder
function. The first non-zero eigenvalue
\begin{equation}
\lambda_{\rm MF}(\eps)=2\beta\,\nu(\xi_0)
\label{eq:lambdaMF}
\end{equation}
is the gap of one independent site. Inserting this $\varphi$ in \eqref{eq:onesiteRQ} and letting
$N\to\infty$ gives
\begin{equation}
\limsup_{N\to\infty}\lambda_1(N,\eps)\ \le\ \lambda_{\rm MF}(\eps),\qquad 1<\eps<2 .
\end{equation}
The function $\lambda_{\rm MF}$ decreases from the sphere value $1/(\eps-1)$ at $\eps\to1$, where
$\varphi$ is linear, to $1/4$ at $\eps\to2^-$. Figure~\ref{fig:lmf_lMC} compares $\lambda_{\rm MF}(\eps)$ with the relaxation rates measured by a Montecarlo numerical algorithm.

We measured the relaxation rates of the Brownian motion directly, with an exactly reversible Montecarlo dynamics: each move rotates three occupations along the circle that preserves their mass and energy, and time is calibrated so that the generator is $\Delta_\M$, the normalisation of $\lambda_{\rm MF}$. Rates were measured separately in the symmetry sectors of $S_N$, extrapolated to zero step size and then to $N\to\infty$; details are in Appendix~\ref{app:numerics}.

The results are reproduced in Figure \ref{fig:lmf_lMC}.

\subsection{\texorpdfstring{The condensate phase $\eps>2$}{The condensate phase eps>2}}

For $\eps>2$ the result of the previous subsection is no longer useful. In the limit $N\to\infty$
the one-site marginal tends to $p(q)=e^{-q}$, the bulk at $\beta=0$, and the condensate
$q_i\simeq a\sqrt N$, $a=\sqrt{\eps-2}$, disappears from it: it carries weight $1/N$ and escapes
to infinity. The solutions of
\begin{equation}
    -\varphi''(q)+\varphi'(q)=\lambda\varphi(q),\quad \varphi'(0)=0,
\end{equation}
that are square integrable with respect to $e^{-q}dq$ are the constant, with $\lambda_0=0$, and a
continuum starting at $\lambda=1/4$. Formally $\lambda_{\rm MF}(\eps>2)=1/4$, which is a correct
but very loose bound.

The slow process in the condensed phase is the transfer of the condensate from one site to
another. Its simplest description uses two sites, say $1$ and $2$, and test functions of $q_1$ and
$q_2$ only. Such a function can be $+1$ when the condensate sits on site $1$, $-1$ when it sits on
site $2$, and $0$ when it sits anywhere else. During the transfer the configuration must cross
the tunnel of Section~\ref{sec:thermo}. This leads to:

\begin{theorem}[upper bound]\label{thm:upper}
Fix $\eps>2$ and let $a=\sqrt{\eps-2}$. Then
\[
\lambda_1(N,\eps)\ \le\ \exp\bigl\{-b(\eps)\sqrt N+O(\ln N)\bigr\},\qquad b(\eps)=(\sqrt2-1)\,a .
\]
\end{theorem}
We prove this theorem by going to two-site test functions. The exponential bound is {\it not} improved by going to 3 and higher variables test functions.

\subsection{Two-site test functions}

The reduction of Section~\ref{sec:fluid} extends to functions of two occupations,
$f(q)=\varphi(q_1,q_2)-\E_\mu\varphi$. The subtraction $-\E_\mu\varphi$ is necessary to orthogonalize with respect to the ground state, constant function. The gradient is
$\nabla_\M f=\sum_{k=1,2}\partial_k\varphi\,\nabla_\M q_k$, and the matrix
$\Pi_{kl}=\langle\nabla_\M q_k,\nabla_\M q_l\rangle$, $k,l\in\{1,2\}$, given in Lemma~\ref{lem:calc}, depends on $q_1$ and $q_2$ only. Both averages in the variational principle therefore reduce
exactly to integrals over the quadrant $(q_1,q_2)$, with the two-site density of
Lemma~\ref{lem:density}(3),
\begin{equation}
p_2(q_1,q_2)=\frac{\Omega_{N-2}\bigl(N-q_1-q_2,\;N\eps-q_1^2-q_2^2\bigr)}{\Omega_N(N,N\eps)} .
\end{equation}
This is the density drawn in Figure~\ref{fig:sections}. One finds
\begin{equation}
\lambda_1\ \le\ \frac{\displaystyle\iint dq_1\,dq_2\;p_2(q_1,q_2)\,\nabla\varphi\cdot \Pi\,\nabla\varphi}
{\displaystyle\iint dq_1\,dq_2\;p_2(q_1,q_2)\,\bigl(\varphi-\E_\mu\varphi\bigr)^2},
\label{eq:twositeRQ}
\end{equation}
for every $\varphi$, at every $N$. In the same way as we replaced $\gamma\to 1$ in the previous subsection, since $\Pi$ is the restriction of an orthogonal projector, $\nabla\varphi\cdot \Pi\,\nabla\varphi\leq |\nabla\varphi|^2$, and $\Pi\to\delta_{kl}$ at the price of a weaker but still valid bound. The loss is $O(1/N)$ where $q_1,q_2=O(1)$, and a factor of order one near the condensate; at the neck, $\Pi$ acts as the identity along the tube direction $(1,-1)$.

The optimal $\varphi$ solves a two-dimensional Sturm--Liouville
problem (with $\Pi$ replaced by $\delta_{kl}$), $-\nabla\cdot(p_2\nabla\varphi)=\Lambda^{(2)}_N\,p_2\,\varphi$ with natural boundary conditions. 

Its first non-zero eigenvalue $\Lambda^{(2)}_N$ is an upper bound on $\lambda_1$ and lies below the one-site value $\Lambda_N$, since one-site functions are a special case.

The two-site problem has a transparent meaning. In the potential $W=-\ln p_2$ the quadrant looks
like panel (b) of Figure~\ref{fig:sections}:
\begin{itemize}
\item two wells along the axes, $q_1\simeq a\sqrt N$ or $q_2\simeq a\sqrt N$, in which the condensate sits on site $1$ or on site $2$;
\item a bulk region near the origin, which represents the condensate sitting on any of the other $N-2$ sites;
\item the tube joining the two wells, whose narrowest point is the neck on the diagonal.
\end{itemize}
The lowest antisymmetric mode is a two-well Kramers problem: $\varphi\simeq+1$ in one well,
$\simeq-1$ in the other, $\simeq0$ in the bulk.

To prove Theorem~\ref{thm:upper} it is enough to write down one explicit function of this type.
We do not solve the two-dimensional problem but a connection with Witten's supersymmetric quantum mechanics \cite{witten1982supersymmetry} is mentioned at the end of this section.

\subsection{The trial function}

Let $q_\pm=a\sqrt{N/2}\pm \ell/2$, with a fixed width $\ell>0$, and let $\chi(q)$ be the ramp equal to
$0$ for $q\le q_-$, to $1$ for $q\ge q_+$, and linear in between. We take
\begin{equation}
\varphi(q_1,q_2)=\chi(q_1)-\chi(q_2).
\label{eq:trial2}
\end{equation}
Then $\varphi=+1$ when site $1$ is above the neck level $a\sqrt{N/2}$ and site $2$ is below it,
$\varphi=-1$ in the opposite case, and $\varphi=0$ when both are below, which is the case when the
condensate is elsewhere. The function changes only when one of the two occupations crosses the neck
level. The three ingredients of the Rayleigh quotient are as follows.

\emph{Mean.} $\varphi$ is antisymmetric under $q_1\leftrightarrow q_2$, which preserves $\mu$. Hence
$\E_\mu\varphi=0$ and $\varphi\perp1$.

\emph{Dirichlet energy.} The intrinsic gradient is the projection of the ambient one, so
$|\nabla_\M\varphi|^2\le \chi'(q_1)^2+\chi'(q_2)^2$. With $I=[q_-,q_+]$ and $p_N$ the one-site
marginal of \eqref{eq:onesiteRQ},
\begin{equation}
\E_\mu|\nabla_\M\varphi|^2\ \le\ \frac{2}{\ell^2}\int_I p_N(q)\,dq\ \le\ \frac2\ell\,\sup_{q\in I}p_N(q).
\label{eq:energy2}
\end{equation}
The energy is thus set by the probability that one given site sits exactly at the neck level.

\emph{Variance.} $\varphi^2=1$ when one of the two sites is above $q_+$ and the other below
$q_-$. Since $\sum_iq_i^2=N\eps$, at most $n_*=N\eps/q_-^2\simeq2\eps/a^2$ sites can
exceed $q_-$. So if a condensate above $q_+$ exists, it can be paired with at least
$N-n_*$ sites below $q_-$. By exchangeability,
\begin{equation}
\E_\mu\varphi^2\ \ge\ 2\,\mu\bigl(q_1\ge q_+,\,q_2\le q_-\bigr)\ \ge\
\frac{2(N-n_*)}{N(N-1)}\;\mu\Bigl(\max_iq_i\ge q_+\Bigr).
\label{eq:var2}
\end{equation}
For $\eps>2$ the last probability tends to one (Lemma~\ref{lem:condensate} below), and
$\E_\mu\varphi^2\ge(2/N)(1-o(1))$.

Inserting \eqref{eq:energy2} and \eqref{eq:var2} in the variational principle gives the two-site
Rayleigh bound
\begin{equation}
\lambda_1\ \le\ \frac N\ell\,\bigl(1+o(1)\bigr)\,\sup_{|q-a\sqrt{N/2}|\le \ell/2}p_N(q).
\label{eq:rayleigh2}
\end{equation}
Everything is now reduced to the value of the one-site marginal at the neck level.

\subsection{The single-site density at the neck level}

The heuristic picture of Section~\ref{sec:thermo} predicts $p_N$ at once. Fix $q_1=z\sqrt N$ with
$z<a$. The remaining sites must carry the excess energy $(a^2-z^2)N$, and the cheapest way to do it
is a second spike of height $\sqrt{(a^2-z^2)N}$. Then
\begin{equation}
p_N(z\sqrt N)\ \approx\ \exp\bigl\{-\sqrt N\,[\Phi_1(z)-a]\bigr\},\qquad
\Phi_1(z)=z+\sqrt{a^2-z^2}.
\end{equation}
$\Phi_1$ is the $\ell^1$ norm of the two spikes at fixed $\ell^2$ norm $a$. It is maximal at
$z=a/\sqrt2$, where the two spikes are equal: the configuration is at the neck of the tube joining
site $1$ to the site of the second spike. There $\Phi_1=\sqrt2\,a$, and
$p_N(a\sqrt{N/2})\approx e^{-b(\eps)\sqrt N}$. Geometrically, the slice $\{q_1=a\sqrt{N/2}\}$
cuts all the $N-1$ tubes leaving the chamber of site $1$, each at its neck. This is the valley
between the bulk and the condensate peak visible in panel (c) of Figure~\ref{fig:sections}.

The rigorous version is:

\begin{proposition}[the neck valley]\label{prop:valley}
Fix $\eps>2$ and $K>0$. Uniformly for $|q-a\sqrt{N/2}|\le K$,
\begin{equation}
p_N(q)\ \le\ \exp\bigl\{-b(\eps)\sqrt N+O(\ln N)\bigr\}.
\label{eq:valley}
\end{equation}
\end{proposition}

\begin{lemma}[a condensate exists]\label{lem:condensate}
Fix $\eps>2$ and $K>0$. Then $\mu\bigl(\max_iq_i<a\sqrt{N/2}+K\bigr)\le\exp\{-b(\eps)\sqrt N+O(\ln N)\}$.
\end{lemma}

The proofs (Appendix~\ref{app:upper}) use the same tools as Proposition~\ref{prop:total}, plus
one observation. Condition on the largest of the occupations $q_2,\dots,q_N$, say $q_2=v$. All the
other sites are then bounded by $v$, and on a bounded support the energy fugacity may be taken
negative: this is the precise sense in which the saddle point sits on the branch cut $\beta<0$ of
\eqref{eq:gsingle}. Lemma~\ref{lem:tilt} then bounds the residual partition function by
$\exp\{N-\sqrt N F(z_1,z_2)+O(\ln N)\}$, with
\begin{equation}
F(z_1,z_2)=z_1+z_2+\frac{(a^2-z_1^2-z_2^2)_+}{\max(z_2,a/2)},\qquad z_i=q_i/\sqrt N ,
\label{eq:Fchernoff}
\end{equation}
with $(x)_+=\max(x,0)$ as in Lemma~\ref{lem:tilt}.

The exponent of \eqref{eq:valley} follows from an elementary inequality: at $z_1=a/\sqrt2$,
\begin{equation}
\min_{z_2\ge0}F\Bigl(\frac a{\sqrt2},z_2\Bigr)=\sqrt2\,a ,
\label{eq:minF}
\end{equation}
attained at $z_2=a/\sqrt2$, i.e.\ at the neck.

\begin{remark}[meaning of the Chernoff exponent]
The term $(a^2-|z|^2)/z_2$ is the cost of storing the excess energy in sites at the maximal allowed
height $z_2$. Inequality \eqref{eq:minF} is then $\|z\|_1\ge\|z\|_2^2/\|z\|_\infty$, applied to
all the sites that carry the excess. It contains the extremisation of $\|z\|_1$ at fixed $\|z\|_2$
of Section~\ref{sec:thermo} and, in addition, rules out all multi-spike configurations at once: a
further spike can never be higher than $z_2$, and splitting the excess among lower sites only
increases $\|z\|_1$. This replaces the qualitative ``concatenation'' argument below
\eqref{eq:budget}. The bound is sharp at the neck at leading order. What is lost is the $O(1)$
factor $e^{-c_*(\eps)}$ of \eqref{eq:tubearea}, which plays no role for the exponent.
\end{remark}

\begin{proof}[Proof of Theorem~\ref{thm:upper}]
Insert \eqref{eq:valley} with $K=\ell/2$ in \eqref{eq:rayleigh2}, for any fixed $\ell$.
\end{proof}

The exponent does not depend on $\ell$: it is a property of the neck. A wider ramp improves only the
prefactor. Expanding $\Phi_1$ around its maximum, $\Phi_1''(a/\sqrt2)=-2\sqrt2/a$, the valley of
$p_N$ is flat within a constant factor over $|q-a\sqrt{N/2}|\lesssim N^{1/4}$, the window already
found in \eqref{eq:neckwidth}. Taking $\ell\simeq N^{1/4}$ gains a factor $N^{-1/4}$.

\begin{remark}[beyond the exponent]\label{rem:kramers}
The optimal trial function changes across each tube rather than along the lines $q_{1,2}=a\sqrt{N/2}$, and its Rayleigh quotient is a capacity. At the tube the residual fluid is in its Gaussian window, $x-2\sim N^{-1/2}\ll N^{-1/3}$, so $\Omega_{N-2}$ follows from Lemma~\ref{lem:llt} with the covariance of $(X,X^2)$, $X\sim\Exp(1)$, whose determinant is $4$. The density across the tube is then a Gaussian of width $1/a$, the integral of $p_2$ across the neck is $e^{-b(\eps)\sqrt N-c_*(\eps)}/(\sqrt2N)$, and the conductance of one tube is
\[
\kappa=\frac{e^{-b(\eps)\sqrt N-c_*(\eps)}}{\sqrt{\sqrt2\,\pi a}\;N^{5/4}} .
\]
The $N$ chambers, of mass $1/N$, form a complete graph with hopping rate $N\kappa$ per pair, whence applying arguments from random processes one can get an estimate (neither an upper nor a lower bound) \cite{Eyring1935,Kramers1940,HanggiTalknerBorkovec1990,BovierEckhoffGayrardKlein2004}
\[
\lambda_1\simeq N^2\kappa=\frac{N^{3/4}}{\sqrt{\sqrt2\,\pi a}}\,e^{-b(\eps)\sqrt N-c_*(\eps)},
\]
which assumes that each chamber equilibrates fast. Functions of one or two sites give twice this value, because $q_1$ crosses the tubes $1$--$j$ at $45^\circ$. For the one-site problem the prediction agrees with the exact evaluation of $p_N$ to about $1\%$ for $400\le N\le6400$ at $\eps=4$.

However, we must now remember that the choice of the Brownian motion (or equivalently of the metric to be put on the microcanonical manifold) has a certain degree of arbitrariness and changing it between various ``natural" choices changes the prefactor of the bound but not the exponent. So one should not take the prefactor seriously when comparing, for example, with Hamiltonian dynamics. 
\end{remark}

\begin{remark}[the exponent]\label{rem:exponent}
Theorems~\ref{thm:lower} and~\ref{thm:upper} give
\[
(\sqrt2-1)\,a\ \le\ \liminf_{N\to\infty}\frac{-\log\lambda_1}{\sqrt N}\ \le\ \limsup_{N\to\infty}\frac{-\log\lambda_1}{\sqrt N}\ \le\ a .
\]
Above $\eps_c$ the spectral gap therefore decays as $e^{-\Theta(\sqrt N)}$: faster than any power of $N$, but slower than $e^{-cN}$ for every $c>0$. This should be compared with the transition in the lowest eigenvalue of the Laplacian on $M$ described in \cite{arezzo2022localization}. The loss in the lower bound does not come from the neck. It comes from the transport of Steps 2--3 of Section~\ref{sec:case2}, which must reach the configurations whose largest coordinate is close to $\eps$ (Lemma~\ref{lem:cover}): the density oscillation along a cap between $q_i\approx2\eps$ and $q_i\approx T$ is $e^{T}$, while the neck costs only $e^{b(\eps)\sqrt N}$. The identification of the exponent, presumably $b(\eps)$ from the upper bound, is left open.
\end{remark}

\emph{Supersymmetric formulation.} During our investigations, we have found that the whole discussion of the upper bound can be recast in the language of Witten's supersymmetric quantum mechanics \cite{witten1982supersymmetry}. The substitution $\psi=\sqrt p\,\varphi$ turns the Sturm--Liouville problems of this section into Schr\"odinger operators $H=A^\dagger A$ with superpotential $W=-\ln p$, whose zero-energy state is the equilibrium. In the condensed phase, with $\hbar=N^{-1/2}$, the exponentially small eigenvalues arise from tunnelling between the wells through the neck, as in the semiclassical analysis of the Witten complex \cite{helffer1985puits,helffer2004quantitative,bovier2005metastability} would follow from boundary and cusp versions of the Eyring--Kramers formula \cite{Berglund2013,BovierDenHollander2015}. We will present this point of view elsewhere.

\section{Lower bounds and a unified Bochner criterion}\label{sec:lower}

The remaining Sections of the paper are devoted to the proof of the following:

\begin{theorem}[lower bounds]\label{thm:lower}\leavevmode
\begin{enumerate}
\item[(i)] If $1<\eps<2$ and $N(2-\eps)^2\ge2000$, then
\(\lambda_1(N,\eps)\ge(2-\eps)^2/[100(\eps-1)]\).
\item[(ii)] There are $c>0$ and $N_0$ such that
\(\lambda_1(N,2)\ge cN^{-4}e^{-16N^{1/3}}\) for $N\ge N_0$.
\item[(iii)] For fixed $\eps>2$ there are $c_\eps,C_\eps>0$ and $N_\eps$ such that, for $N\ge N_\eps$,
\begin{equation}\label{eq:lower-supercritical}
\lambda_1(N,\eps)\ge c_\eps N^{-C_\eps}
\exp\!\left\{-\sqrt{N(\eps-2)}-16\left(\frac{N}{\eps-2}\right)^{1/4}\right\}.
\end{equation}
\end{enumerate}
\end{theorem}
The statements (ii)--(iii) are asymptotic, because the caps $G_i$ must be nonempty: at $\eps=2$ the condition $T<t^*$, i.e.\ $16N^{1/3}<1+\sqrt{N-1}$, already requires $N$ of order $1.7\cdot10^7$.

Part (i) follows from a Bochner--Reilly criterion (Proposition~\ref{prop:reilly}) applied to $M$ with a single vector field (Section~\ref{sec:global}). For $\eps\ge2$ the criterion is applied to the caps $G_i$, and the resulting Poincar\'e inequalities are transported to the sets $U_i$ and glued (Section~\ref{sec:route}).

\subsection{The Bochner--Reilly criterion}\label{sec:criterion}

We will use the following version of Reilly's application of the (integrated) Bochner's Formula for domain with corners.

\begin{proposition}[Bochner--Reilly criterion with corners]\label{prop:reilly}\label{prop:singular}
Let $D$ be a compact connected $d$-\hspace{0pt}dimensional Riemannian domain with smooth corners, $d\ge3$, and let $F\subset\overline D$ be finite (possibly empty). Suppose that $X$ is smooth on $\overline D\setminus F$, that $|X|\le C\,\operatorname{dist}(\cdot,F)^{-1}$ near $F$, and, for some $\lambda_*>0$,
\begin{align}
 \operatorname{Ric}+(\operatorname{div}X-|X|^2)g&\ge\lambda_*g
 &&\text{on }D\setminus F,\label{eq:criterion-interior}\\
 \mathrm{II}&\ge\langle X,\nu\rangle g
 &&\text{on each open boundary face outside }F.\label{eq:criterion-boundary}
\end{align}
Here $\nu$ is outward, $\mathrm{II}(v,v)=\langle\nabla_v\nu,v\rangle$, and \eqref{eq:criterion-boundary} is an inequality on vectors tangent to the face. Then $\lambda_1(D)\ge\lambda_*$. If $F=\varnothing$, the restriction $d\ge3$ is unnecessary.
\end{proposition}

\begin{proof}
Integrals over $D$ are taken with respect to its Riemannian measure $\sigma=\cH^d$, and integrals over a boundary face with respect to $\cH^{d-1}$. Let $u$ be a first nonconstant variational Neumann eigenfunction, write $\lambda_1=\lambda_1(D)$, and set $V=\nabla u$. The regularity theorem of Mitrea--Taylor--Vasy \cite[Proposition~3.2 and Corollary~4.2]{MitreaTaylorVasy2005}, applied to the absolute one-form $du$, gives $V\in H^1(D)$: indeed $d(du)=0$, $\delta du=\lambda_1u\in L^2$, and $\nu\mathbin{\lrcorner}du=0$. The facewise Gaffney identity (the integrated Bochner--Weitzenb\"ock formula; see \cite[Theorem~3]{RaulotSavo2011} for one-forms with zero normal component) gives
\begin{equation}\label{eq:gaffney}
 \lambda_1\int_D|V|^2\,d\sigma
 =\int_D\bigl(|\nabla V|^2+\operatorname{Ric}(V,V)\bigr)\,d\sigma
   +\sum_S\int_S\mathrm{II}(V,V)\,d\cH^{d-1},
\end{equation}
where the sum is over the open boundary faces $S$. For a smooth field tangent to every face, the same identity with left-hand side $\int_D((\operatorname{div}V)^2+|dV^\flat|^2)\,d\sigma$ follows by applying the divergence theorem to $(\operatorname{div}V)V-\nabla_VV$. This identity extends to $H^1$ tangential fields by density and the trace theorem: in a corner chart, extend each coordinate component oddly across its own face and evenly across the other faces, then mollify. Substitution of $V=\nabla u$ yields \eqref{eq:gaffney}.

Integrating $|\nabla V+Z\otimes V|^2\ge0$, for any Lipschitz field $Z$, gives
\begin{equation}\label{eq:square}
 \int_D|\nabla V|^2\,d\sigma\ge
 \int_D(\operatorname{div}Z-|Z|^2)|V|^2\,d\sigma
 -\sum_S\int_S\langle Z,\nu\rangle|V|^2\,d\cH^{d-1}.
\end{equation}
If $F=\varnothing$, take $Z=X$ in \eqref{eq:square} and combine \eqref{eq:gaffney} with \eqref{eq:criterion-interior}--\eqref{eq:criterion-boundary}. Otherwise put $r=\operatorname{dist}(\cdot,F)$ and $Z=\zeta_\delta X$, where $\zeta_\delta=0$ on $r<\delta$, $\zeta_\delta=1$ on $r>2\delta$, $0\le\zeta_\delta\le1$, and $|\nabla\zeta_\delta|\le C/\delta$. Since Ricci curvature and the face second fundamental forms are bounded, \eqref{eq:criterion-interior}--\eqref{eq:square} imply
\begin{equation}\label{eq:cutoff}
\begin{split}
 \lambda_1\int_D|V|^2\,d\sigma\ge{}&\lambda_*\int_D|V|^2\,d\sigma
 -C\int_{D\cap\{r<2\delta\}}|V|^2\,d\sigma\\
 &-C\int_{D\cap\{\delta<r<2\delta\}}r^{-2}|V|^2\,d\sigma
 -C\sum_S\int_{S\cap\{r<2\delta\}}|V|^2\,d\cH^{d-1}.
\end{split}
\end{equation}
The unweighted volume and face error integrals in \eqref{eq:cutoff} tend to zero by integrability and the trace theorem. For the integral with weight $r^{-2}$, use the local point-singularity Hardy inequality: for $p\in\overline D$ and sufficiently small fixed $a>0$,
\begin{equation}\label{eq:hardy}
 \int_{D\cap B_a(p)}\frac{|V|^2}{\operatorname{dist}(q,p)^2}\,d\sigma(q)
 \le C_a\|V\|_{H^1(D)}^2.
\end{equation}
Estimate \eqref{eq:hardy} follows from the Euclidean Hardy inequality \cite[Section~1.2.3]{BEL}, after a cutoff and a bounded $H^1$ extension in a corner chart \cite[Chapter~VI, Section~3]{Stein1970}. It requires $d\ge3$ but no vanishing boundary trace; the $H^1$ norm includes the necessary $L^2$ term. Since $F$ is finite, \eqref{eq:hardy} makes $r^{-2}|V|^2$ locally integrable near $F$. The weighted integral in \eqref{eq:cutoff} therefore tends to zero as $\delta\downarrow0$, proving $\lambda_1\ge\lambda_*$.
\end{proof}

\emph{On the microcanonical surface.} Proposition~\ref{prop:reilly} is applied to $D=M$ for $1<\eps<2$, and to $D=G_i$ for $\eps\ge2$. By Lemma~\ref{lem:calc}(2)--(3) its two hypotheses become scalar conditions:
\begin{align}
\frac{N-3}{R^2}+\dv X-|X|^2&\ge\lambda_* ,\label{eq:micro-interior}\\
\langle X,\nabla q_j\rangle&\ge\frac{1-K}{R^2}
\quad\text{on a face }\{q_j=K\}\text{ of }\{q_j\ge K\}.\label{eq:micro-boundary}
\end{align}
On the faces $\{q_j=0\}$ ($K=0$), \eqref{eq:micro-boundary} reads $\langle X,\nabla q_j\rangle\ge1/R^2$: since $\nabla q_j$ points into the domain, $X$ must enter through each of these faces with normal component at least $1/(R^2|\nabla q_j|)$, which compensates their concavity. On the convex face $\{q_i=T\}$ of $G_i$ ($K=T>1$), \eqref{eq:micro-boundary} holds for every field tangent to the face.

\emph{The cost of the concavity.} Integrating \eqref{eq:micro-interior} over $M$, with the divergence theorem, \eqref{eq:micro-boundary} on the faces $\{q_j=0\}$, and the coarea identity $\int_{\{q_j=0\}}|\nabla q_j|^{-1}\,d\cH^{N-3}=\sigma(M)\,p_N(0)$, where $p_N$ is the one-site marginal of \eqref{eq:onesiteRQ}, gives for every field allowed in Proposition~\ref{prop:reilly}
\begin{equation}\label{eq:budget-bound}
\lambda_*\ \le\ \frac{N-3}{R^2}-\frac{N\,p_N(0)}{R^2}-\int_M|X|^2\,d\mu\ \le\ \frac{1-p_N(0)}{\eps-1}.
\end{equation}
As $N\to\infty$, $p_N(0)$ tends to the density at $0$ of the one-site law of Section~\ref{sec:thermo}: this limit is smaller than $1$ for $\eps<2$, with $1-\lim_Np_N(0)=(2-\eps)/2+O((2-\eps)^2)$, and equals $1$ for $\eps\ge2$. A single field on $M$ can therefore give an $N$-uniform bound only for $\eps<2$. For $\eps\ge2$ the criterion is applied to the caps $G_i$, whose slices are below the critical energy (Remark~\ref{rem:margin}).

\subsection{Corners}\label{sec:cornersub}
Proposition~\ref{prop:reilly} is stated for compact domains with smooth corners, i.e. points with a neighbourhood diffeomorphic to a neighbourhood of $0$ in $[0,\infty)^m\times\R^{d-m}$. The following observation proves that smooth corners can fail only where the faces through a point have linearly dependent gradients and that this happens only at the points of $F_\eps:=\{q\in M:\ q_j\in\{0,\eps\}\ \forall j\},$
the configurations with $k$ coordinates equal to $\eps$ and the others equal to $0$, which exist only at the exceptional energies $N/\eps=k\in\Z$, and on the face $\{q_i=K\}$ of a cap for finitely many levels $K$.

\begin{proposition}[corners of the domains]\label{prop:corners}
\begin{enumerate}[label=\textup{(\alph*)}]
\item Every $q\in M\setminus F_\eps$ has a neighbourhood in $M$ diffeomorphic to a neighbourhood of $0$ in $[0,\infty)^m\times\R^{d-m}$, $m$ being the number of coordinates vanishing at $q$. In particular, if $N/\eps\notin\Z$, $M$ is a compact connected domain with smooth corners.
\item For $K\in(\eps,t^*)$ outside a set of at most $2N$ values, the cap $\{q\in M:\ q_i\ge K\}$ is a compact connected domain with smooth corners, and it contains no point of $F_\eps$. This holds for $G_i$, by the choice of $T$ in \eqref{eq:threshold}.
\item If $N/\eps=k\in\Z$ with $k\le N-2$, $M$ is not a domain with corners at the points $p\in F_\eps$: on $M$
\begin{equation}\label{eq:cusp}
\sum_{j:\,p_j=0}q_j=\frac{|q-p|^2}{2\eps},
\end{equation}
so every curve in $M$ starting at $p$ is tangent to the $(k-1)$-dimensional subspace $\{v\in T_p:\ v_j=0\text{ if }p_j=0\}$, and $k-1<d$. \textup{(}For $k=N-1$, $M=S_R$ is a round sphere.\textup{)}
\end{enumerate}
\end{proposition}
The proof is in Appendix~\ref{app:corners}.

Part (b) also covers the slices of a cap at exceptional energies: at a level $t$ with $n/e(t)\in\Z$, the faces of the slice $\Sigma_i(t)\simeq M_{n,e(t)}$ meet non-transversally at the points of $c(t)F_{e(t)}$, but these points are not in $F_\eps$, since $q_i=t>\eps$, so they are regular corner points of the cap. In Step~1 they are only isolated singular points of the vector field, which Proposition~\ref{prop:reilly} allows. At the exceptional energies $N/\eps\in\Z$ of $M$ itself, part (c) shows that Proposition~\ref{prop:reilly} cannot be applied to $M$; Theorem~\ref{thm:lower}(i) is obtained there by approximating the energy from below (Section~\ref{sec:completion}).

\subsection{\texorpdfstring{The global field for $1<\eps<2$}{The global field for 1<eps<2}}\label{sec:global}
The vector field $X$ is built from the profile
\begin{equation}\label{eq:profile-definition}
\varphi(x)=\frac{(2-\eps)\,x\,(x-\eps)^2}{(2-\eps)x+8\eps},\qquad x\ge0,
\end{equation}
through the averages
\begin{equation}\label{eq:field-definition}
H=\avg{\varphi},\qquad
J=(\eps-1)\varphi'(0)+\avg{(x-\eps)\varphi'},\qquad
X=\frac{\nabla H}{J}.
\end{equation}
Since $\varphi\ge0$ vanishes only at $0$ and $\eps$, $H>0$ on $M\setminus F_\eps$.

\begin{proposition}[the global field]\label{prop:global}
Let $1<\eps<2$ and $N(2-\eps)^2\ge2000$. The field $X=\nabla H/J$ of \eqref{eq:field-definition} is smooth on $M\setminus F_\eps$, and there
\begin{align}
\langle X,\nabla q_j\rangle&=\frac1{R^2}\qquad\text{on each face }\{q_j=0\},\label{eq:global-flux}\\
\frac{N-3}{R^2}+\dv X-|X|^2&\ge\frac{(2-\eps)^2}{100(\eps-1)},\label{eq:global-curvature}\\
|X|&\le2\,(NH)^{-1/2}.\label{eq:global-size}
\end{align}
In particular, if $F_\eps=\varnothing$, Proposition~\ref{prop:reilly} applies with $D=M$, $F=\varnothing$ and $\lambda_*=(2-\eps)^2/[100(\eps-1)]$, and gives $\lambda_1(N,\eps)\ge(2-\eps)^2/[100(\eps-1)]$.
\end{proposition}

Thus \eqref{eq:global-flux} is \eqref{eq:micro-boundary} with equality, \eqref{eq:global-curvature} is \eqref{eq:micro-interior}, and \eqref{eq:global-size} controls $X$ near the zeros of $H$.

To prove the Proposition we will need the following elementary properties of the function $\varphi$, which will be proved in Appendix~\ref{app:profile}. For $x\ge0$:
\begin{enumerate}[label=\textup{($\varphi$\arabic*)},ref=$\varphi$\arabic*]
\item\label{phi:basic} $0\le\varphi(x)\le(x-\eps)^2$, $\varphi'(0)=\eps(2-\eps)/8$ and $|\varphi''(x)|\le2$;
\item\label{phi:J} $\displaystyle 2\varphi(x)\le(x-\eps)\Big(\varphi'(x)-\varphi'(0)\Big(1-\frac x\eps\Big)\Big)\le\frac{13}4\,\varphi(x)$;
\item\label{phi:gradH} $\displaystyle\Big(\varphi'(x)-\varphi'(0)\Big(1-\frac x\eps\Big)\Big)^2\le\frac{13}4\,(x-\eps)\Big(\varphi'(x)-\varphi'(0)\Big(1-\frac x\eps\Big)\Big)$;
\item\label{phi:gradJ} $\displaystyle\Big(\varphi'(x)+(x-\eps)\varphi''(x)-\big(\varphi'(0)-\eps\varphi''(0)\big)\Big(1-\frac x\eps\Big)\Big)^2\le144\,\varphi(x)$;
\item\label{phi:moment} for every probability measure $d\nu$ on $[0,\infty)$ with $\int x\,d\nu=1$ and $\int x^2\,d\nu=\eps$,
\[
\varphi'(0)+\int(\varphi''-\varphi')\,d\nu\ \ge\ \frac{3(2-\eps)^2}{16\eps}.
\]
\end{enumerate}

\begin{proof}
All statements are on $M\setminus F_\eps$, where $H>0$.

\emph{The normalisation.} Since $\avg{x}=1$ and $\avg{x^2}=\eps$, we have $\avg{(x-\eps)(1-x/\eps)}=1-\eps$, so that $J=\avg{(x-\eps)\big(\varphi'(x)-\varphi'(0)(1-x/\eps)\big)}$. By (\ref{phi:J}), and by $H\le\avg{(x-\eps)^2}=\eps(\eps-1)<2(\eps-1)$, which follows from (\ref{phi:basic}),
\begin{equation}\label{eq:J-estimate}
2H\le J\le\tfrac{13}4H<\tfrac{13}2(\eps-1).
\end{equation}
In particular $J>0$, so $X$ is smooth on $M\setminus F_\eps$.

\emph{The flux \eqref{eq:global-flux}.} On $\{q_j=0\}$, Lemma~\ref{lem:calc}(4) with $\psi=\varphi$ gives $\langle\nabla H,\nabla q_j\rangle=\frac1N\bigl(\varphi'(0)+\avg{(x-\eps)\varphi'}/(\eps-1)\bigr)=J/R^2$, because $R^2=N(\eps-1)$. Dividing by $J$ gives \eqref{eq:global-flux}.

\emph{Gradients and the size \eqref{eq:global-size}.} Apply Lemma~\ref{lem:calc}(4) to $H=\avg{\varphi}$ with $\alpha+\beta x=\varphi'(0)(1-x/\eps)$, and to $J=\avg{\psi}$, $\psi(x)=(x-\eps)\big(\varphi'(x)-\varphi'(0)(1-x/\eps)\big)$, with $\alpha+\beta x=\psi'(0)(1-x/\eps)$. Since $\psi'(x)-\psi'(0)(1-x/\eps)=\varphi'(x)+(x-\eps)\varphi''(x)-\big(\varphi'(0)-\eps\varphi''(0)\big)(1-x/\eps)$, properties (\ref{phi:gradH}) and (\ref{phi:gradJ}) give, using $H\le J/2$,
\begin{equation}\label{eq:gradient-estimates}
|\nabla H|^2\le\frac{13}{4}\,\frac JN,\qquad
|\nabla J|^2\le144\,\frac HN\le72\,\frac JN.
\end{equation}
With $J\ge2H$, the first bound gives $|X|^2=|\nabla H|^2/J^2\le13/(8NH)<4/(NH)$, which is \eqref{eq:global-size}.

\emph{The interior inequality \eqref{eq:global-curvature}.} Lemma~\ref{lem:calc}(5) gives
$J+(\eps-1)\Delta H=(\eps-1)\{\varphi'(0)+\avg{\varphi''-\varphi'}\}+\mathcal R_\varphi$, where $|\mathcal R_\varphi|\le8(\eps-1)/N$ by (\ref{phi:basic}). Applying (\ref{phi:moment}) to the measure $d\nu(x)=\frac1N\sum_j\delta(x-q_j)\,dx$, whose first two moments are $1$ and $\eps$, and using $\eps<2$,
\begin{equation}\label{eq:laplacian-estimate}
J+(\eps-1)\Delta H\ge(\eps-1)\left(\frac{3(2-\eps)^2}{32}-\frac8N\right).
\end{equation}
Since $X=\nabla H/J$ and $R^2=N(\eps-1)$,
\begin{equation}\label{eq:curvature-expansion}
\frac{N-3}{R^2}+\dv X-|X|^2
=\frac{J+(\eps-1)\Delta H}{(\eps-1)J}
 -\frac3{N(\eps-1)}
 -\frac{\langle\nabla H,\nabla J\rangle+|\nabla H|^2}{J^2}.
\end{equation}
By \eqref{eq:J-estimate}, $1/(\eps-1)<13/(2J)$; by \eqref{eq:gradient-estimates}, $|\nabla H|\,|\nabla J|\le\sqrt{(13/4)\cdot72}\,J/N<16J/N$, since $(13/4)\cdot72=234<256$. Hence the last two terms in \eqref{eq:curvature-expansion} satisfy
\begin{equation}\label{eq:curvature-errors}
\frac3{N(\eps-1)}<\frac{20}{NJ},\qquad
\frac{\langle\nabla H,\nabla J\rangle+|\nabla H|^2}{J^2}<\frac{16+4}{NJ}=\frac{20}{NJ}.
\end{equation}
Inserting \eqref{eq:laplacian-estimate} and \eqref{eq:curvature-errors} in \eqref{eq:curvature-expansion},
\[
\frac{N-3}{R^2}+\dv X-|X|^2\ \ge\ \frac1J\left(\frac{3(2-\eps)^2}{32}-\frac{48}{N}\right).
\]
Since $N(2-\eps)^2\ge2000$, $48/N<50/N\le(2-\eps)^2/40$, and $3/32-1/40=11/160>13/200$. With $J<\frac{13}2(\eps-1)$ this gives
\[
\frac{N-3}{R^2}+\dv X-|X|^2\ \ge\ \frac{13}{200}\,\frac{(2-\eps)^2}J\ \ge\ \frac{(2-\eps)^2}{100(\eps-1)}.
\]

\emph{The last claim.} If $F_\eps=\varnothing$, $X$ is smooth on $M$, and $M$ is a domain with smooth corners by Proposition~\ref{prop:corners}(a); Proposition~\ref{prop:reilly} applies with $F=\varnothing$.
\end{proof}

\begin{figure}[tb]
\centering
\includegraphics[width=0.92\columnwidth]{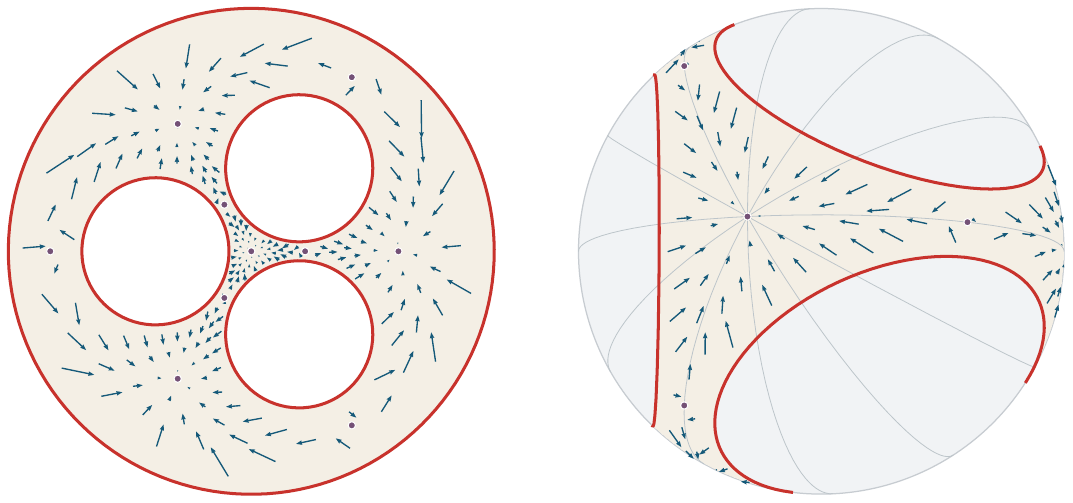}
\caption{The global field $X=\nabla H/J$ of Proposition~\ref{prop:global} on the two-dimensional surface $M_{4,1.7}$. The field crosses the concave faces $\{q_j=0\}$ inwards, with the flux \eqref{eq:global-flux} that compensates their concavity.}
\label{fig:global}
\end{figure}

Proposition~\ref{prop:global} proves Theorem~\ref{thm:lower}(i) when $N/\eps\notin\mathbb Z$; Section~\ref{sec:completion} handles the exceptional energies. Figure~\ref{fig:global} shows the field $X$ on the two-dimensional surface $M_{4,1.7}$.

\subsection{The route from local estimates to the global gap}\label{sec:route}
Recall the cells and the caps of Section~\ref{sec:regions},
\begin{gather*}
A_i=\{q\in\Mo:\ q_i>q_j\ \forall j\ne i\},\qquad
U_i=\{q\in\Mo:\ q_i>\eps\},\\
H_i=\{q\in\Mo:\ q_i>2\eps\},\qquad
G_i=\{q\in M:\ q_i\ge T\},
\end{gather*}
where $T$, given by \eqref{eq:threshold}, lies slightly above $\sqrt{N(\eps-2)}$ for $\eps>2$ and equals $16N^{1/3}$ up to $1$ for $\eps=2$; thus $A_i\subset U_i$ and $G_i\subset H_i\subset U_i$. For $\eps\ge2$ the proof bounds the Poincar\'e constants $P_G=P(G_i)$, $P_H=P(H_i)$, $P_U=P(U_i)$ in four steps:
\begin{enumerate}
\item[Step 1.] Lift the vector field $X$ in \eqref{eq:field-definition} to the slices of $G_i$. The cap estimate \eqref{eq:cap-poincare} gives $P_G\le C_\eps N^{1/2}$ for $\eps>2$ and $P_G\le CN^{2/3}$ for $\eps=2$.
\item[Step 2.] Replace $f$ by the function $g$ in \eqref{eq:frozen-function}, which is constant for $q_i\ge T+1$. The convex chart of $H_i$ then gives $P_H\le C_\eps N^3e^T$ in \eqref{eq:H-final}.
\item[Step 3.] Integrate along meridians across $\eps<q_i\le2\eps$ to obtain $P_U\le C_\eps(1+P_H)$ in \eqref{eq:U-final}.
\item[Step 4.] Apply Lemma~\ref{lem:glue} to the cells $A_i\subset U_i$ and the overlap bound $\mu(U_i\cap U_j)\ge\alpha$. Equation~\eqref{eq:global-gap} gives $\lambda_1\ge c_\eps\alpha N^{-4}e^{-T}$.
\end{enumerate}
Only Step~1 uses Proposition~\ref{prop:reilly}. The factor $e^T$ enters through the density ratio bounded in \eqref{eq:density-ratio-bound}; Proposition~\ref{prop:overlap} bounds $\alpha^{-1}$ by a power of $N$.

\section{\texorpdfstring{From deep caps to the covering sets, $\eps\ge2$}{From deep caps to the covering sets, eps>=2}}\label{sec:case2}
Constants in this section may depend on the fixed energy $\eps$, but not on $N$, except for the local singularity estimates explicitly identified in Section~\ref{sec:step1}. Domains carry the induced metric and the area measure $\sigma$, and $T$ is the level \eqref{eq:threshold}.

\subsection{\texorpdfstring{Step 1: the lifted field on $G_i$}{Step 1: the lifted field on G\_i}}\label{sec:step1}
On $M_{n,e(t)}$, write $y_j=q_j/c(t)$ for $j\ne i$ and set
\begin{equation}\label{eq:slice-field}
\begin{gathered}
\varphi_t(y)=\frac{\theta(t)y(y-e(t))^2}{\theta(t)y+8e(t)},\qquad
H_t=\frac1n\sum_{j\ne i}\varphi_t(y_j),\\
J_t=(e(t)-1)\varphi_t'(0)+\frac1n\sum_{j\ne i}(y_j-e(t))\varphi_t'(y_j),\qquad
Z_t=\frac{\nabla_{M_{n,e(t)}}H_t}{J_t}.
\end{gathered}
\end{equation}
Thus $Z_t$ is the vector field $X$ in \eqref{eq:field-definition} with $(N,\eps)$ replaced by $(n,e(t))$, so that $2-\eps$ becomes $\theta(t)$. If $D_t:M_{n,e(t)}\to\Sigma_i(t)$ is the homothety of ratio $c(t)$, define
\begin{equation}\label{eq:lifted-field}
Y=c^{-2}(D_t)_*Z_t,\qquad X=\varpi Y,\qquad \varpi=\frac{\rho^2}{cR^2}\in(0,1].
\end{equation}
The bound $\varpi\le1$ follows from $cR^2-\rho^2=N(t-1)(t-\eps)/(N-1)\ge0$.
By Remark~\ref{rem:margin}, Proposition~\ref{prop:global} applies on $M_{n,e(t)}$ with $(N,\eps)$ replaced by $(n,e(t))$, away from the finite set $F_{e(t)}$. Scaling \eqref{eq:global-curvature} by $D_t$ gives
\begin{equation}\label{eq:scaled-curvature}
\frac{n-3}{\rho^2}+\dv_\Sigma Y-|Y|^2\ge\frac{n\theta^2}{100\rho^2}.
\end{equation}
Since $X$ is tangent to the slices and $\varpi$ depends only on $t$, the volume element in Lemma~\ref{lem:slices}(2) gives $\dv X=\varpi\dv_\Sigma Y$. Combining \eqref{eq:lifted-field}--\eqref{eq:scaled-curvature} yields
\begin{equation}\label{eq:cap-curvature}
\begin{aligned}
\frac{N-3}{R^2}+\dv X-|X|^2
&\ge\frac1{R^2}\left(\frac{n\theta^2}{100c}-D(t)\right)
\ge\frac{n\theta^2}{200R^2}
\ge\frac{\kappa_N}{400(\eps-1)},\\
D(t)&=\frac{N-4}{c}-(N-3)
=\frac{(t-2)N-3t+4}{N-t}\le2t.
\end{aligned}
\end{equation}
The second inequality in \eqref{eq:cap-curvature} uses $c\le1$ and $n\theta^2\ge400t$ from \eqref{eq:slice-margin}; the third uses $n\ge N/2$ and $R^2=N(\eps-1)$. On each face $\{q_j=0\}$, $j\ne i$, scaling \eqref{eq:global-flux} gives
\begin{equation}\label{eq:cap-flux}
\langle X,\nabla q_j\rangle=\varpi\langle Y,\nabla_\Sigma q_j\rangle
=\frac{\varpi}{c\,n(e-1)}=\frac1{R^2}.
\end{equation}
On $q_i=T$, tangency gives $\langle X,\nabla q_i\rangle=0\ge(1-T)/R^2$. Thus \eqref{eq:cap-flux} and tangency verify \eqref{eq:micro-boundary} on every face of $G_i$.

To check the singularity hypothesis of Proposition~\ref{prop:reilly}, let $F$ consist of the pole $P_i$ and the configurations at the finitely many levels $t_k\in[T,t^*)$ with $n/e(t_k)=k\in\mathbb Z$, $k<n$, at which $k$ of the coordinates $y_j$, $j\ne i$, equal $e(t_k)$ and the others vanish. Near such a configuration, let $A\subset\{1,\ldots,N\}\setminus\{i\}$ be the set of its positive coordinates. The function $H_t$ in \eqref{eq:slice-field} satisfies
\begin{equation}\label{eq:slice-H-lower}
H_t\ge c_0\left(\sum_{j\notin A,\,j\ne i}y_j+
\sum_{j\in A}(y_j-e(t))^2\right).
\end{equation}
The mass constraint gives
$k(e(t)-e(t_k))=-\sum_{j\notin A,\,j\ne i}y_j-\sum_{j\in A}(y_j-e(t))$.
Since $e'(t_k)\ne0$, this identity and \eqref{eq:slice-H-lower} imply $H_t\ge c_1\dist(q,F)^2$.
Estimate \eqref{eq:global-size}, with the slice parameters, gives $|Z_t|\le2(nH_t)^{-1/2}$; hence \eqref{eq:lifted-field} implies
\begin{equation}\label{eq:cap-singularity}
|X(q)|\le C/\dist(q,F).
\end{equation}
Near $P_i$, all $y_j\to1$, so \eqref{eq:slice-field} gives $H_t\ge c_2e(t)(e(t)-1)$ and $|X|\le C\rho/(cR^2)\to0$. The constants in this paragraph may depend on $N$: Proposition~\ref{prop:reilly} requires the inverse-distance bound for each fixed $N$.

The cap $G_i$ is a compact connected domain with smooth corners (Proposition~\ref{prop:corners}(b)). Proposition~\ref{prop:reilly}, applied using \eqref{eq:cap-curvature}, \eqref{eq:cap-flux}, and \eqref{eq:cap-singularity}, gives
\begin{equation}\label{eq:cap-poincare}
\lambda_1(G_i)\ge\frac{\kappa_N}{400(\eps-1)},\qquad
P_G\le\begin{cases}
(\eps-1)\sqrt N/\sqrt{\eps-2},&\eps>2,\\
N^{2/3}/90,&\eps=2.
\end{cases}
\end{equation}

\begin{figure}[tb]
\centering
\includegraphics[width=0.92\columnwidth]{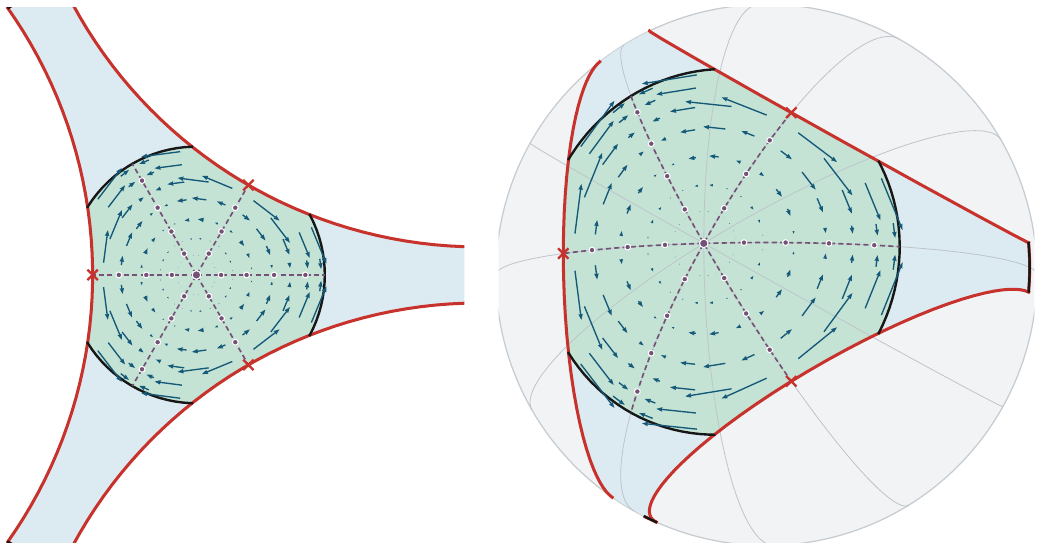}
\caption{The lifted field of Step~1 on an invariant two-sphere of $M_{10^8,2.5}$, where the cap $G_1$ is nonempty. The corner of $G_1$, the singular points of the field and the pole are marked.}
\label{fig:local}
\end{figure}

Figure~\ref{fig:local} shows the lifted field on an invariant two-sphere of $M_{10^8,2.5}$, where the cap $G_1$ is nonempty, and marks its corner, its singular points and the pole.

\subsection{\texorpdfstring{Step 2: freezing and transport from $G_i$ to $H_i$}{Step 2: freezing and transport from G\_i to H\_i}}\label{sec:step2}
For $K>\eps$ the coordinates $q_j$, $j\ne i$, measured as fractions of their total mass, give a convex chart of the cap $\{q_i>K\}$. We state it for $H_i$.

\begin{lemma}\label{lem:chart}
Fix $i$ and write $t:=q_i$. On $H_i$, $t$ ranges over $(2\eps,t^*]$, the value $t^*$ being attained only at the pole $P_i$. Let $\mathcal T=\{x\in[0,\infty)^{N-1}:\ \sum_jx_j=1\}$ and $k(t)=(N\eps-t^2)/(N-t)^2$; $k$ is strictly decreasing on $(\eps,N)$. The weighted coordinates $x_j=q_j/(N-t)$, $j\ne i$, map $H_i$ diffeomorphically onto the convex set
\[
\Omega=\{x\in\mathcal T^\circ:\ |x|^2<k(2\eps)\},
\]
a simplex intersected with a ball. On $\Omega$, $t$ is the function $t(x)=k^{-1}(|x|^2)$, and the inverse chart is $F(x)=\big(q_i=t(x),\ q_j=(N-t(x))\,x_j\big)$. Moreover, for $N\ge7\eps$:
\begin{enumerate}[label=\textup{(\alph*)}]
\item $d\sigma=w(t(x))\,dx$, where $dx$ is the $(N-2)$-dimensional Lebesgue measure on the hyperplane $\{\sum_jx_j=1\}$ and
\[
w(t)=\sqrt{\tfrac{\eps-1}{N-1}}\ \frac{(N-t)^{N-1}}{t-\eps};
\]
$w$ is decreasing in $t$, and $N^{N-1}e^{-t-t^2/N}\le(N-t)^{N-1}\le N^{N-1}e^{-t+t/N}$ for $0\le t\le N/2$;
\item $\diam(\Omega)^2\le8\eps/N$;
\item $\|dF\|^2\le(2+6/\eps)\,N^3$.
\end{enumerate}
\end{lemma}
The main factor $(N-t)^{N-1}\approx N^{N-1}e^{-t}$ of the density says that taking mass $t$ away from the other $N-1$ coordinates shrinks their volume by $(1-t/N)^{N-1}\approx e^{-t}$; it is the origin of the factor $e^T$ below. The proof is in Appendix~\ref{app:chart}.

Let $m$ be the mean of $f\in H^1(H_i)$ on $G_i$. Choose $0\le\chi\le1$ with $\chi=0$ on $(-\infty,T]$, $\chi=1$ on $[T+1,\infty)$ and $|\chi'|\le2$, and set
\begin{equation}\label{eq:frozen-function}
g=f-\chi(q_i)(f-m).
\end{equation}
Since $|\nabla q_i|\le1$, the Poincar\'e bound \eqref{eq:cap-poincare} gives
\begin{equation}\label{eq:freezing-estimates}
\begin{aligned}
\int_{H_i}|f-g|^2d\sigma&\le P_G\int_{H_i}|\nabla f|^2d\sigma,\\
\int_{H_i}|\nabla g|^2d\sigma&\le(2+8P_G)\int_{H_i}|\nabla f|^2d\sigma.
\end{aligned}
\end{equation}
In the convex chart $F:\Omega\to H_i$ of Lemma~\ref{lem:chart}, $g\circ F=m$ for $t\ge T+1$. For the density $w(t)$ of Lemma~\ref{lem:chart}(a), define
\begin{equation}\label{eq:capped-density}
\widetilde w(t)=\max\{w(t),w(T+1)\},\qquad
\frac{\sup_\Omega\widetilde w}{\inf_\Omega\widetilde w}\le
\Xi:=\frac{w(2\eps)}{w(T+1)}
=\left(\frac{N-2\eps}{N-T-1}\right)^{N-1}\frac{T+1-\eps}{\eps}.
\end{equation}
The Payne--Weinberger inequality \cite{PW,Beb} on $\Omega$ and \eqref{eq:capped-density} give the weight $\widetilde w$ a Poincar\'e constant at most $\Xi\operatorname{diam}(\Omega)^2/\pi^2$. The same constant bounds the variance of $g\circ F$ with weight $w$: use $w\le\widetilde w$ in the variance integral and $\nabla(g\circ F)=0$ wherever $w\ne\widetilde w$ in the energy integral. Lemma~\ref{lem:chart}(b)--(c) and \eqref{eq:freezing-estimates} therefore give
\begin{equation}\label{eq:H-preliminary}
\begin{aligned}
\min_{c\in\R}\int_{H_i}|g-c|^2d\sigma&\le C_\eps N^2\Xi\int_{H_i}|\nabla g|^2d\sigma,\\
P_H&\le2P_G+C_\eps N^2\Xi(1+P_G).
\end{aligned}
\end{equation}
The logarithmic bounds in Lemma~\ref{lem:chart}(a), applied to \eqref{eq:capped-density}, imply
\begin{equation}\label{eq:density-ratio-bound}
\Xi\le\frac{T+1}{\eps}
\exp\!\left(T+1-2\eps+\frac{2\eps+(T+1)^2}{N}\right)
\le\begin{cases}
C_\eps N^{1/2}e^T,&\eps>2,\\
CN^{1/3}e^T,&\eps=2.
\end{cases}
\end{equation}
Substituting \eqref{eq:cap-poincare} and \eqref{eq:density-ratio-bound} into \eqref{eq:H-preliminary} proves
\begin{equation}\label{eq:H-final}
P_H\le C_\eps N^3e^T.
\end{equation}
\subsection{\texorpdfstring{Step 3: meridian transport from $H_i$ to $U_i$}{Step 3: meridian transport from H\_i to U\_i}}\label{sec:step3}
Write $d\sigma=\upsilon(t)\,dt\,d\omega$ on $U_i$, as in \eqref{eq:meridian-density}.
By Lemma~\ref{lem:slices}(3), the admissible meridian directions increase with $t$. Thus every $(t,\omega)\in L:=U_i\setminus H_i$ is joined within $U_i$ to $(s,\omega)$ for $2\eps<s<2\eps+1$, and
\begin{equation}\label{eq:meridian-difference}
|f(t,\omega)-f(s,\omega)|^2
\le2(\eps+1)\int_t^s|\nabla f(r,\omega)|^2\,dr.
\end{equation}
Let $m$ be the mean of $f$ on $H_i$. Average
$|f(t,\omega)-m|^2\le2|f(s,\omega)-m|^2+2|f(t,\omega)-f(s,\omega)|^2$
over $s\in(2\eps,2\eps+1)$ and integrate over $L$. Equations~\eqref{eq:meridian-difference} and~\eqref{eq:meridian-density} give
\begin{equation}\label{eq:layer-estimate}
\int_L|f-m|^2d\sigma
\le C_\eps\int_{H_i}|f-m|^2d\sigma
+C_\eps\int_{U_i}|\nabla f|^2d\sigma.
\end{equation}
Adding $\int_{H_i}|f-m|^2d\sigma$ to \eqref{eq:layer-estimate} and applying \eqref{eq:H-final} on $H_i$ proves
\begin{equation}\label{eq:U-final}
P_U\le C_\eps(1+P_H)\le C_\eps N^3e^T.
\end{equation}
For $f\in H^1(U_i)$, equations~\eqref{eq:meridian-difference}--\eqref{eq:layer-estimate} follow by local Sobolev approximation and absolute continuity along almost every meridian.

\section{Completion of the lower bounds}\label{sec:completion}\label{sec:corners}
\noindent\textbf{Exceptional energies in Theorem~\ref{thm:lower}(i).}
Suppose $\eps_0=N/k\in(1,2)$ and $N(2-\eps_0)^2\ge2000$. On the unit sphere $S$ in $\{\sum_jv_j=0\}$ put $K_\tau=\{u\in S:u_j\ge-\tau\ \forall j\}$. The map $u\mapsto\mathbf1+R(\eps)u$ identifies $M_{N,\eps}$ with $(K_{1/R(\eps)},R(\eps)^2g_S)$. As nonexceptional $\eps\uparrow\eps_0$, the sets $K_{1/R(\eps)}$ decrease to $K_{1/R(\eps_0)}$, and their difference has vanishing measure. For restrictions of smooth functions on $S$, the Poincar\'e inequality obtained from \eqref{eq:global-flux} and \eqref{eq:global-curvature} therefore passes to the limit. Its spectral-gap lower bound in the metric $g_S$ is
\begin{equation}\label{eq:limit-gap}
R(\eps)^2\frac{(2-\eps)^2}{100(\eps-1)}
=\frac{N(2-\eps)^2}{100}.
\end{equation}
Restrictions of smooth functions on $S$ are dense in $H^1(K_{1/R(\eps_0)})$: truncate the function, cut it off near the finite set corresponding to $F_{\eps_0}$ at squared-gradient cost $O(\delta^{d-2})$, where $d=N-2>2$, and use corner-chart density away from that set. The limit of \eqref{eq:limit-gap}, rescaled by $R(\eps_0)^{-2}$, proves Theorem~\ref{thm:lower}(i) at $\eps_0$. This limiting argument treats the points of $F_{\eps_0}$, at which $M$ is not a domain with corners (Proposition~\ref{prop:corners}(c)).

\noindent\textbf{Step 4: Theorem~\ref{thm:lower}(ii)--(iii).}
The local Poincar\'e inequalities on the sets $U_i$ are glued by the following lemma.

\begin{lemma}\label{lem:glue}
Let $A_i\subset U_i$, $1\le i\le N$, where the $A_i$ partition $M$ up to null sets and $\mu(A_i)=1/N$. If $P(U_i)\le P$ and $\mu(U_i\cap U_j)\ge\alpha>0$ for $i\ne j$, then
\begin{equation}\label{eq:gluing}
\Var_\mu f\le N\bigl(1+\sqrt{2/\alpha}\bigr)^2P\int_M|\nabla f|^2\,d\mu.
\end{equation}
\end{lemma}
\begin{proof}
Let $m_i$ be the mean of $f$ on $U_i$, $v=m_i$ on $A_i$, and $E_i=\int_{U_i}|\nabla f|^2d\mu$. The Poincar\'e inequalities on $U_i$ imply
\begin{equation}\label{eq:gluing-local}
\|f-v\|_2^2\le P\sum_iE_i\le NP\mathcal E(f,f),\qquad
\alpha|m_i-m_j|^2\le2P(E_i+E_j).
\end{equation}
Summing the second inequality in \eqref{eq:gluing-local} gives
$\Var_\mu v=N^{-2}\sum_{i<j}(m_i-m_j)^2\le2P\mathcal E(f,f)/\alpha$.
The triangle inequality for the centred $L^2$ norm, together with the first inequality in \eqref{eq:gluing-local}, gives $(\sqrt N+\sqrt{2/\alpha})^2P$ as a Poincar\'e constant, which is no larger than the constant in \eqref{eq:gluing}.
\end{proof}

Lemma~\ref{lem:cover} supplies the cover $U_i$ and the cells $A_i$ required by Lemma~\ref{lem:glue}. Proposition~\ref{prop:overlap} gives $\alpha\ge\alpha_2>0$ at $\eps=2$ and $\alpha\ge N^{-C_\eps}$ at fixed $\eps>2$. Since $\alpha\le1$ and $(1+\sqrt{2/\alpha})^2\le8/\alpha$, substituting \eqref{eq:U-final} into \eqref{eq:gluing} yields
\begin{equation}\label{eq:global-gap}
\begin{aligned}
\Var_\mu f&\le C_\eps N^4\alpha^{-1}e^T\int_M|\nabla f|^2d\mu,\\
\lambda_1(N,\eps)&\ge c_\eps\alpha N^{-4}e^{-T}.
\end{aligned}
\end{equation}
Finally, \eqref{eq:threshold} gives $e^{-T}\ge e^{-1}e^{-T_0}$. Substitution of the two values of $T_0$ and the corresponding overlap bounds into \eqref{eq:global-gap} proves Theorem~\ref{thm:lower}(ii)--(iii).
\FloatBarrier


\section*{Concluding summary}

We have shown that substituting the Hamiltonian dynamics with a Brownian dynamics on the microcanonical surface, we can put rigorous bounds on the equilibration time. The fluid phase, for $\eps<2$ equilibrates in $O(1)$ time. The condensate phase for $\eps>2$ equilibrates in time which diverges with $N$,
\begin{equation}
    \teq\asymp e^{b_*(\eps)\sqrt{N}+o(\sqrt N)}
\end{equation}
and we have the exact bounds
\begin{equation}
    (\sqrt{2}-1)\sqrt{\eps-2}\leq b_*(\eps)\leq\sqrt{\eps-2}.
\end{equation}

To obtain the lower bound on $\lambda_1$, which is an upper bound to $\teq$, we have used a geometric technique which we dubbed {\it Bochner method} because it rests on Bochner's formula: integrated for the gradient of an eigenfunction, it turns a lower bound on the Ricci curvature into a lower bound on $\lambda_1$, and a vector field that crosses the concave faces of $M$ inwards turns the wrong-sign boundary term into an interior cost that the curvature of the sphere can pay (Proposition~\ref{prop:reilly}). For $\eps<2$ one such field on the whole surface gives a bound uniform in $N$; for $\eps\ge2$ it works only on the deep caps around a condensate, and transporting and gluing these local inequalities across $M$ gives $\lambda_1\ge e^{-\sqrt{N(\eps-2)}-O(N^{1/4})}$ (Theorem~\ref{thm:lower}). The proof of the lower bound to $\lambda_1$ is involved (although conceptually simple) but it is necessary so. In fact, one must prove that {\it for any square integrable, zero average, function $f$} one cannot get a ratio $\int d\mu |\nabla_M f|^2/\int d\mu |f|^2$ lower than something prescribed (the putative bound), which is exponentially small in $\sqrt{N}$. Physically, this means that {\it the only} long time dynamics occurs because of the entropic effects at the neck. In other words, that the chambers equilibrate quickly. This is the most difficult part of the proof.

To obtain the upper bound to $\lambda_1$, which is a lower bound on $\teq$, instead it is sufficient to present an appropriate variational trial functions on the microcanonical surface. We found that for $\eps\leq2$ it is sufficient (as supported by the numerics) to have trial functions of a single variable $q_1$ while for $\eps>2$ one must switch to two functions $(q_1,q_2)$ in order to capture the event of a localized mass to shift from one lattice site to another. When we cannot solve the 2-$d$ partial differential equation, we can present a simple linear trial function of $(q_1,q_2)$ which does the job.

Our upper bound also admits a natural formulation in terms of Witten's supersymmetric quantum mechanics \cite{witten1982supersymmetry}: the equilibrium measure is the unique zero-energy state, and in the condensed phase the $N-1$ exponentially small eigenvalues are the ground states of the $N$ chambers, lifted from zero by tunnelling through the necks. This point of view will be developed in a separate work.

The choice of the Brownian motion also deserves a comment. Any reversible dynamics with the microcanonical measure as equilibrium has the same geometric bottleneck, so the exponent $b(\eps)\sqrt N$ in the condensed phase is a property of the measure and not of the dynamics. The prefactors and the $O(1)$ corrections, instead, depend on the metric used to define the Brownian motion, and for $\eps\le2$ the dependence affects the gap itself. For instance, the Brownian motion in the variables $\psi_i$, on the surface $\{\sum_i|\psi_i|^2=N,\ \sum_i|\psi_i|^4=N\eps\}\subset\mathbb C^N$, has the same equilibrium measure but a diffusion coefficient proportional to $q_i$ in the occupations, closer to the Hamiltonian mass transfer between sites. At $\eps=2$ this surface is, to leading order, an Einstein manifold with $\mathrm{Ric}\simeq2$, and its gap is close to $2$, whereas for the Brownian motion in the variables $q_i$ studied here $\lambda_1\lesssim1/4$ (with the equality most likely being correct, as shown by the numerics). 

The $\sqrt{N}$ dependence of the (log of) the equilibration time is an {\it entropy effect} as measured by the reduction of the volume of the region (the {\it tunnel}) that the phase space point has to pass to go from one condensate configuration to another one. Such entropic barriers are common to other statistical mechanics problems \cite{schweizer2003entropic,reguera2006entropic}, and in recent years have been argued to hinder the performance of both classical and quantum algorithms \cite{bellitti2021entropic}. 

Going back to look at the Hamiltonian dynamics, we must notice a curious ``coincidence" with the work \cite{iubini2019PRL}. In \cite{iubini2019PRL} the Hamiltonian dynamics of the DNLSE coupled to a bath, at positive but high temperature has been studied and a breather has been found to melt on timescales $\tau_{b}\sim e^{\alpha_{\rm melt} b(0)}$ where $b(0)=|\psi_i(0)|^2$ is the mass of the breather at position $i$ and $\alpha_{\rm melt}=0.91$. If we bootstrap this numerical law to our configurations at negative temperatures, where $\eps>2$, $b(0)=q_i(0)=\sqrt{\eps-2}\sqrt{N}$ and $\ln\tau_{\rm eq}=\alpha\sqrt{\eps-2}\sqrt{N}$ with $\alpha\in[\sqrt{2}-1,1]$ in this paper so we have the same $\sqrt{N}$ dependence and a compatible prefactor. The comparison is meaningful because the background of \cite{iubini2019PRL} has mean mass density close to one, so that $b(0)$ and $q_i$ are measured in the same units, and their exponent $\alpha_{\rm melt}$ is in the interval we found. The mechanisms are {\it prima facie} different: in \cite{iubini2019PRL} the slowness comes from an adiabatic invariant associated with the fast rotation of the breather, which our Brownian approximation discards, while our barrier is purely entropic and has no counterpart at positive temperature. However, the results suggest they must be connected somehow.

The slowness of the Hamiltonian dynamics above $\eps_c=2$ is also visible in the numerical study \cite{mithun2018weakly}, where the infinite-temperature line $\eps=2$ of the Gross--Pitaevskii lattice was found to mark the onset of long-lived localized excitations, whose lifetimes have power-law distributions and make the dynamics weakly nonergodic at slightly larger energy densities.

To our knowledge, the one presented in this paper is the first rigorous result in which the relaxation time of a dynamics on the DNLS microcanonical surface is bounded, from above and from below, by quantities that diverge with the system size: in the condensed phase both bounds are stretched exponentials, $e^{\Theta(\sqrt N)}$. This should be contrasted with the rigorous Hamiltonian results based on adiabatic invariants \cite{carati2012exponentially,giorgilli2015extensive}, which are uniform in $N$ and in which the small parameter is the specific energy, and with the time scales of condensate motion in stochastic particle systems \cite{BeltranLandim}, which grow polynomially in the system size. 

We believe the techniques exposed in this paper can be used to find rigorous bounds of this kind, for many stochastic problems in the future.

\section*{Acknowledgements}
Claude Opus was used in this paper heavily to run the numerics, to set up and bring to closure the cumbersome calculations needed in the proof of the lower bound, and for checking overall consistency of the methods and results. We like to believe, however, that the ideas which have not appeared in the literature before are our own.

\appendix

\section{The tilted bound (Lemma~\ref{lem:tilt})}\label{app:tilt}
\begin{proof}[Proof of Lemma~\ref{lem:tilt}]
Write $s=\zeta\sqrt N\ge\frac12$, $S=N-\varsigma\sqrt N$ and $E=N(\eps-r^2)$. For real $\beta$ let $p_{s,\beta}(q)=e^{-q-\beta q^2}\mathbf 1_{[0,s]}(q)/\cZ_s(1,\beta)$, with $\cZ_s(1,\beta)=\int_0^s e^{-q-\beta q^2}dq$. Since $\prod_i e^{-q_i-\beta q_i^2}=e^{-S-\beta E}$ on the support of the $\delta$-functions,
\begin{equation}
\Omega^{<s}_n(S,E)=e^{\,S+\beta E}\,\cZ_s(1,\beta)^n\,f^{s,\beta}_n(S,E)
\label{eq:tilt}
\end{equation}
exactly, where $f^{s,\beta}_n$ is the density of $(\sum_i q_i,\sum_i q_i^2)$ for $n$ i.i.d.\ variables with density $p_{s,\beta}$. By Lemma~\ref{lem:density}(1),
\begin{equation}
\sup f^{s,\beta}_n\ \le\ \frac{\pi}{\sqrt3}\,\bigl(\sup_q p_{s,\beta}(q)\bigr)^3 .
\label{eq:f3}
\end{equation}
\emph{Case $r\ge a$.} Take $\beta=0$. Then $\cZ_s(1,0)\le1$, and $\sup p_{s,0}\le1/(1-e^{-1/2})$ since $s\ge\frac12$; so $\ln\Omega^{<s}_n\le S+O(1)=N-\varsigma\sqrt N+O(1)$.

\emph{Case $r<a$.} Take $\beta=-\kappa/\sqrt N<0$ with
\[
\kappa=\frac{1-\delta}{c},\qquad c=\max\Bigl(\zeta,\frac a2\Bigr),\qquad \delta=\frac{6\ln N}{a\sqrt N}.
\]
The exponent $\mathcal L(q)=-q+\kappa q^2/\sqrt N$ is convex, and at the cutoff $\mathcal L(s)=-s(1-\kappa\zeta)\le0$ because $\kappa\zeta\le1-\delta$. So $e^{\mathcal L}\le1$ on $[0,s]$, $\cZ_s(1,\beta)\ge1-e^{-1/2}$, and \eqref{eq:f3} is bounded. Split $\cZ_s(1,\beta)=\int_0^se^{\mathcal L(q)}\,dq$ at $L=N^{1/8}$:
\begin{itemize}[nosep]
\item on $[0,L]$, $e^{\kappa q^2/\sqrt N}=1+\kappa q^2/\sqrt N+O(q^4/N)$ gives $\int_0^Le^{\mathcal L(q)}\,dq\le1+2\kappa/\sqrt N+O(N^{-1})$;
\item on $[L,s]$, convexity bounds $\mathcal L$ by its endpoint values: $\mathcal L(L)\le-L/2$, while $\mathcal L(s)\le-s/2$ if $\zeta\le a/4$ (where $\kappa\zeta\le\frac12$) and $\mathcal L(s)\le-\delta s\le-\frac32\ln N$ if $\zeta>a/4$; so this part is $O(N^{-1})$.
\end{itemize}
Hence $n\ln\cZ_s(1,\beta)\le2\kappa\sqrt N+O(1)$. Inserting $S$ and $\beta E=-\kappa\sqrt N(\eps-r^2)$ in \eqref{eq:tilt}, and using $\eps-2=a^2$,
\[
\ln\Omega^{<s}_n\ \le\ N-\sqrt N\big[\varsigma+\kappa(a^2-r^2)\big]+O(1)\ =\ N-\sqrt N\Big[\varsigma+\frac{a^2-r^2}{c}\Big]+\frac{\delta\sqrt N(a^2-r^2)}{c}+O(1),
\]
and the middle term is at most $2a\delta\sqrt N=12\ln N$.
\end{proof}

\section{\texorpdfstring{Proofs of Propositions~\ref{prop:total} and~\ref{prop:overlap}}{Proofs of Propositions III.8 and III.9}}\label{app:Ceps}
\begin{proof}[Proof of Proposition~\ref{prop:total}]
\emph{Lower bound.} Integrating out one coordinate in the definition of $\Omega_N$ gives $\Omega_N(N,N\eps)=\int_0^\infty\Omega_{N-1}(N-x,N\eps-x^2)\,dx$; restrict $x$ to $[a\sqrt N,a\sqrt N+1]$. The other $n=N-1$ coordinates have mass $S'=N-x$ and energy $E'=N\eps-x^2$, so that $S'-n=-a\sqrt N+O(1)$ and $E'-2n=O(\sqrt N)$. Lemma~\ref{lem:llt} gives $f_n(S',E')\ge c/N$, hence $\Omega_n(S',E')=e^{S'}f_n(S',E')\ge c\,e^{N-a\sqrt N}/N$. Integrating over $x$,
\[
\Omega_N(N,N\eps)\ \ge\ c\,e^{N-a\sqrt N}/N .
\]

\emph{Upper bound.} The largest coordinate exceeds $\eps$ (Lemma~\ref{lem:cover}) and is almost surely attained at a single site. Conditioning on its position and value $u$,
\[
\Omega_N(N,N\eps)=N\int_\eps^{t^*}\Omega^{<u}_{N-1}\big(N-u,\,N\eps-u^2\big)\,du .
\]
Apply Lemma~\ref{lem:tilt} with $k=1$ and $\varsigma=r=\zeta=u/\sqrt N$, which is admissible since $\eps\le u\le t^*\le\sqrt{N\eps}$ for $N$ large. Its proof (Appendix~\ref{app:tilt}) gives in fact the sharper bound
\[
\ln\Omega^{<u}_{N-1}\big(N-u,\,N\eps-u^2\big)\ \le\ N-\sqrt N\Big[\zeta+(1-\delta)\,\frac{(a^2-\zeta^2)_+}{\max(\zeta,a/2)}\Big]+O(1),\qquad \delta=\frac{6\ln N}{a\sqrt N}.
\]
This exponent lies below $N-a\sqrt N$ by at least:
\begin{itemize}
\item $u-a\sqrt N$ for $\zeta\ge a$;
\item $(1-2\delta)(a\sqrt N-u)$ for $a/2\le\zeta<a$, because $\zeta+(1-\delta)(a^2/\zeta-\zeta)-a=\delta(\zeta-a)+(1-\delta)\,a(a-\zeta)/\zeta\ge(1-2\delta)(a-\zeta)$;
\item $a\sqrt N/4$ for $\zeta<a/2$, because there the bracket is concave in $\zeta$ and at least $\min\big(2(1-\delta)a,\ a/2+\tfrac32(1-\delta)a\big)\ge\tfrac54a$ for $N$ large.
\end{itemize}
Integrating over $u$ gives the stronger upper bound
\[
\Omega_N(N,N\eps)\ \le\ C'_\eps N\,e^{N-a\sqrt N} .
\]
The two bounds give \eqref{eq:normlower}, and the bounds on $f_N(N,N\eps)=e^{-N}\Omega_N(N,N\eps)$ follow.
\end{proof}

\begin{proof}[Proof of Proposition~\ref{prop:overlap}]
By Lemma~\ref{lem:density}(3),
\[
\Omega_N(N,N\eps)\,\mu(U_i\cap U_j)=\int_\eps^\infty\!\!\int_\eps^\infty\Omega_{N-2}\big(N-u-v,\,N\eps-u^2-v^2\big)\,du\,dv .
\]

(a) Let $\eps=2$ and restrict $u,v$ to $[3,4]\subset(\eps,\infty)$. Since $\Omega_n=e^Sf_n$,
\[
\mu(U_i\cap U_j)\ \ge\ \frac1{f_N(N,2N)}\int_3^4\!\!\int_3^4e^{-u-v}f_{N-2}\big(N-u-v,\,2N-u^2-v^2\big)\,du\,dv .
\]
Both densities are evaluated within $O(1)$ of their barycentres $(n,2n)$, where Lemma~\ref{lem:llt} gives $nf_n\to G(0)>0$, $G$ the limiting Gaussian density. Hence the right side tends to $(e^{-3}-e^{-4})^2=2\alpha_2$, and it exceeds $\alpha_2$ for $N\ge N_0$.

(b) Let $\eps>2$ and $a=\sqrt{\eps-2}$. Restrict $u,v$ to $(\eps,\eps+1)$, and write $\Omega_{N-2}(S,E)=\int_0^\infty\Omega_{N-3}(S-x,E-x^2)\,dx$, restricting $x$ to $[a\sqrt N,a\sqrt N+1]$. The remaining $n=N-3$ sites have mass $S'=N-u-v-x$ and energy $E'=N\eps-u^2-v^2-x^2$, so that $S'-n=-a\sqrt N+O(1)$ and $E'-2n=O(\sqrt N)$. Lemma~\ref{lem:llt} gives $f_n(S',E')\ge c/N$, hence $\Omega_n(S',E')=e^{S'}f_n(S',E')\ge c\,e^{N-a\sqrt N}/N$, and
\[
\Omega_N(N,N\eps)\,\mu(U_i\cap U_j)\ \ge\ c\,e^{N-a\sqrt N}/N .
\]
Together with the stronger upper bound $\Omega_N(N,N\eps)\le C'_\eps N\,e^{N-a\sqrt N}$ proved above,
\[
\mu(U_i\cap U_j)\ \ge\ c_\eps N^{-2},
\]
so (b) holds with any $C_\eps>2$.
\end{proof}

Consequently Theorem~\ref{thm:lower}(iii) holds with $N^{-6}$ in place of $N^{-C_\eps}$. We do not know whether the power $N^{-2}$ is sharp; at $\eps=2$ the overlap is bounded below by a constant.

\section{Proofs of the neck estimates}\label{app:upper}

\begin{proof}[Proof of Proposition~\ref{prop:valley}]
\emph{Step 1: exact representation.} Almost surely the largest of $q_2,\dots,q_N$ is attained at a
unique site. By permutation symmetry and \eqref{eq:muconstraint}, conditioning on its position and
value $v$ gives
\begin{equation}
p_N(u)=\frac{N-1}{\Omega_N(N,N\eps)}\int dv\;
\Omega^{<v}_{N-2}\bigl(N-u-v,\;N\eps-u^2-v^2\bigr),
\label{eq:pNexact}
\end{equation}
with $\Omega^{<v}_{N-2}$ as in Lemma~\ref{lem:density}. The range of
$v$ is bounded. On one side $v\le\sqrt{N\eps}$. On the other side the $N-2$ remaining sites carry
mass $N-u-v$ and are all below $v$, so $v\ge(N-u)/(N-1)\ge\frac12$ for $u\le\sqrt{N\eps}$ and $N$
large. Note that $v$ may be larger or smaller than $u$: no ordering between sites $1$ and $2$ is
assumed.

\emph{Step 2: the tilted bound.} Write $z_1=u/\sqrt N$, $z_2=v/\sqrt N$. Lemma~\ref{lem:tilt}
with $k=2$, $\varsigma=z_1+z_2$, $r^2=z_1^2+z_2^2$ and cutoff $\zeta=z_2\ge1/(2\sqrt N)$ gives
\begin{equation}
\ln\Omega^{<v}_{N-2}\bigl(N-u-v,\,N\eps-u^2-v^2\bigr)\ \le\ N-\sqrt N\,F(z_1,z_2)+C\ln N ,
\end{equation}
with $F$ as in \eqref{eq:Fchernoff}.

\emph{Step 3: the variational inequality.} Let $z_1=a/\sqrt2+\delta$ with
$|\delta|\le K/\sqrt N$. Since $\max(z_2,a/2)\ge a/2$, one has
$|\partial F/\partial z_1|\le1+4z_1/a$ on both sides of the circle $|z|=a$, and $F$ is continuous
across it. Hence $F(z_1,z_2)\ge F(a/\sqrt2,z_2)-C_1K/\sqrt N$, with $C_1=1+4\sqrt\eps/a$. At
$z_1=a/\sqrt2$ there are three cases.
\begin{itemize}
\item $z_2\ge a/\sqrt2$: the point is outside the disk, and $F=a/\sqrt2+z_2\ge\sqrt2\,a$.
\item $a/2\le z_2<a/\sqrt2$: $F=a/\sqrt2+z_2+(a^2/2-z_2^2)/z_2=a/\sqrt2+a^2/(2z_2)>\sqrt2\,a$.
\item $0\le z_2<a/2$: $F=a/\sqrt2+a+z_2-2z_2^2/a$. This is concave in $z_2$, so it is at least its
value at the endpoints, $a/\sqrt2+a>\sqrt2\,a$.
\end{itemize}
Hence $\min_{z_2}F(a/\sqrt2,z_2)=\sqrt2\,a$, attained at the neck $z_2=a/\sqrt2$, which is
\eqref{eq:minF}. The integrand of \eqref{eq:pNexact} is therefore at most
$\exp\{N-\sqrt2\,a\sqrt N+C_1K+C\ln N\}$.

\emph{Step 4: normalisation.} By Proposition~\ref{prop:total},
$\Omega_N(N,N\eps)\ge\exp\{N-a\sqrt N-O(\ln N)\}$. The prefactor $N-1$ and the length
$\sqrt{N\eps}$ of the range of $v$ are polynomial. Hence
$p_N(u)\le\exp\{-(\sqrt2-1)a\sqrt N+O(\ln N)\}$, uniformly for $|u-a\sqrt{N/2}|\le K$.
\end{proof}

\begin{proof}[Proof of Lemma~\ref{lem:condensate}]
Let $u_K=a\sqrt{N/2}+K$. As in the proof of the upper bound in Proposition~\ref{prop:total},
condition on the position and value $u$ of the largest coordinate, which exceeds $\eps$ by
Lemma~\ref{lem:cover}:
\[
\Omega_N(N,N\eps)\,\mu\Bigl(\max_iq_i<u_K\Bigr)=N\int_\eps^{u_K}
\Omega^{<u}_{N-1}\bigl(N-u,\,N\eps-u^2\bigr)\,du .
\]
Lemma~\ref{lem:tilt} with $k=1$ and $\varsigma=r=\zeta=u/\sqrt N$ bounds the integrand by
$\exp\{N-\sqrt N\,F(\zeta)+C\ln N\}$, with $F(\zeta)=\zeta+(a^2-\zeta^2)_+/\max(\zeta,a/2)$. On
the range $\zeta\le a/\sqrt2+K/\sqrt N$ there are two cases.
\begin{itemize}
\item $a/2\le\zeta$: $F=a^2/\zeta\ge\sqrt2\,a-O(K/\sqrt N)$.
\item $\zeta<a/2$: $F-a=(a-\zeta)(a+2\zeta)/a$, which is concave and equal to $a$ at both ends,
so $F\ge2a$.
\end{itemize}
Hence the integrand is at most $\exp\{N-\sqrt2\,a\sqrt N+O(\ln N)\}$. Dividing by the lower bound
of Proposition~\ref{prop:total} proves the lemma.
\end{proof}

\section{\texorpdfstring{The corners of the domains (Proposition~\ref{prop:corners})}{The corners of the domains}}\label{app:corners}
\begin{proof}[Proof of Proposition~\ref{prop:corners}]
\emph{Linear independence of the gradients.} The tangent space of $S_R$ at $q$ is the orthogonal complement of $\one$ and $q$, and $\nabla q_j$ is the orthogonal projection of $e_j$ onto it. Hence, for a set $I$ of indices, a relation $\sum_{j\in I}c_j\nabla q_j=0$ holds if and only if $\sum_{j\in I}c_je_j=\alpha\one+\beta q$ for some $\alpha,\beta\in\R$, that is,
\begin{equation}\label{eq:relation}
c_j=\alpha+\beta q_j\quad(j\in I),\qquad 0=\alpha+\beta q_l\quad(l\notin I).
\end{equation}
In the two cases below some $l\notin I$ exists. If $\beta=0$, then \eqref{eq:relation} gives $\alpha=0$ and all $c_j=0$. A nontrivial relation therefore has $\beta\ne0$, and then all the coordinates $q_l$, $l\notin I$, are equal to the same $m=-\alpha/\beta$, which is positive.
\begin{enumerate}[label=\textup{(\roman*)}]
\item Let $q\in M\setminus F_\eps$ and $I=\{j:\ q_j=0\}$; some $l\notin I$ exists because $\sum_jq_j=N$. If the $k$ coordinates outside $I$ were all equal to $m$, the constraints would give $km=N$ and $km^2=N\eps$, so $m=\eps$ and $q\in F_\eps$. Hence the gradients of the coordinates vanishing at $q$ are linearly independent.
\item Let $K\in(\eps,t^*)$, $q\in M$ with $q_i=K$, and $I=\{i\}\cup\{j:\ q_j=0\}$; some $l\notin I$ exists, since otherwise $q_i=N$ and $N^2=N\eps$. If the $k\in\{1,\dots,N-1\}$ coordinates outside $I$ were all equal to $m$, the constraints would give $K+km=N$ and $K^2+km^2=N\eps$, that is, $kK^2+(N-K)^2=kN\eps$. For each $k$ this quadratic equation has at most two roots $K$. Hence, for $K$ outside a set of at most $2N$ values, at every $q\in M$ with $q_i=K$ the gradients of $q_i$ and of the coordinates vanishing at $q$ are linearly independent.
\end{enumerate}

\emph{(a) and (b).} Let $D$ be $M$ or the cap, written as $D=\{q\in S_R:\ f_1(q)\ge0,\dots,f_s(q)\ge0\}$, where the $f$'s are the coordinates and, for the cap, $q_i-K$. $D$ is closed in the compact manifold $S_R$, hence compact; it is connected by Lemma~\ref{lem:cover} and Lemma~\ref{lem:slices}(3). Let $q\in D$, $q\notin F_\eps$, and let $f_1,\dots,f_m$ be the functions vanishing at $q$. Their differentials at $q$ are linearly independent: by (i) if $q_i>K$ or $D=M$, by (ii) if $q_i=K$. (The cap contains no point of $F_\eps$, whose coordinates are at most $\eps<K$.) By the inverse function theorem, $f_1,\dots,f_m$ together with $d-m$ suitable linear functions are coordinates of $S_R$ near $q$. The other $f$'s are positive near $q$, so in these coordinates $D$ is, near $q$, the set $\{f_1\ge0,\dots,f_m\ge0\}$: a neighbourhood of $0$ in $[0,\infty)^m\times\R^{d-m}$.

\emph{(c).} Let $Z=\{j:\ p_j=0\}$ and write $q=p+v$ with $q\in M$. Then $\sum_jv_j=0$, and $\sum_j(p_j+v_j)^2=N\eps=\sum_jp_j^2$ gives $2\eps\sum_{j\notin Z}v_j+|v|^2=0$. Since $\sum_{j\notin Z}v_j=-\sum_{j\in Z}v_j$ and $v_j=q_j$ for $j\in Z$, this is \eqref{eq:cusp}. Along a curve $q(s)$ in $M$ with $q(0)=p$, the left side of \eqref{eq:cusp} is a sum of nonnegative terms of order $s^2$, so $\dot q_j(0)=0$ for $j\in Z$. The vectors of $T_p$ with $v_j=0$ for $j\in Z$ form a subspace of dimension $k-1<d$, while the tangent vectors at $0$ of curves in $[0,\infty)^m\times\R^{d-m}$ span $\R^d$.
\end{proof}

\section{\texorpdfstring{The properties of the profile $\varphi$ (Section~\ref{sec:global})}{The properties of the profile}}\label{app:profile}
\begin{proof}[Proof of (\ref{phi:basic})--(\ref{phi:moment})]
Put
\[
u(x)=\frac{\eps(10-\eps)}{(2-\eps)x+8\eps},\qquad 0<u(x)\le u(0)=\frac{10-\eps}8<\frac98 ,
\]
and note that $\varphi(x)/(x-\eps)=(2-\eps)x(x-\eps)/\big((2-\eps)x+8\eps\big)$ is a rational function, with
\[
\Big(\frac{\varphi(x)}{x-\eps}\Big)^2=\varphi(x)\,\frac{(2-\eps)x}{(2-\eps)x+8\eps}\le\varphi(x).
\]
Let $\xi(x)=(x-\eps)\big(\varphi'(x)-\varphi'(0)(1-x/\eps)\big)$; then $\xi'(x)-\xi'(0)(1-x/\eps)$ is the function in (\ref{phi:gradJ}). Differentiating \eqref{eq:profile-definition} gives
\begin{align*}
&\varphi'(0)=\frac{\eps(2-\eps)}8,\qquad \varphi''=2-\frac{16u^3}{10-\eps},\\
&\varphi'(x)-\varphi'(0)\Big(1-\frac x\eps\Big)=\frac{\varphi(x)}{x-\eps}\cdot\frac{18-\eps+8u(x)}8,\qquad\text{hence}\qquad \xi=\varphi\,\frac{18-\eps+8u}8,\\
&\xi'(x)-\xi'(0)\Big(1-\frac x\eps\Big)=\frac{\varphi(x)}{x-\eps}\,B(x),\qquad B(x)=\frac{\eps^2-24\eps+172}{32}+\frac{(14-\eps)\,u(x)}4+2u(x)^2 .
\end{align*}
The bounds follow from $1<\eps<2$ and $0<8u\le10-\eps<9$.

\emph{(\ref{phi:basic}).} Since $0\le(2-\eps)x/\big((2-\eps)x+8\eps\big)<1$, $0\le\varphi\le(x-\eps)^2$. Moreover $2>\varphi''\ge2-(10-\eps)^2/32>-1$.

\emph{(\ref{phi:J}).} Since $16<18-\eps+8u<26$, $2\varphi\le\xi\le\frac{13}4\varphi$.

\emph{(\ref{phi:gradH}).} By the bound on $\big(\varphi/(x-\eps)\big)^2$,
\[
\Big(\varphi'(x)-\varphi'(0)\Big(1-\frac x\eps\Big)\Big)^2\le\varphi\,\Big(\frac{18-\eps+8u}8\Big)^2=\xi\,\frac{18-\eps+8u}8\le\frac{13}4\,\xi .
\]

\emph{(\ref{phi:gradJ}).} The three terms of $B$ are positive and smaller than $5$, $4$ and $3$, so $0<B<12$ and $\big(\xi'(x)-\xi'(0)(1-x/\eps)\big)^2\le144\,\varphi(x)$.

\emph{(\ref{phi:moment}).} Let $\Psi=\varphi'(0)+\varphi''-\varphi'$. Since $\varphi''-2$ is a negative multiple of $u^3$, whose derivatives alternate in sign, $\varphi^{(4)}<0<\varphi^{(5)}$ and $\Psi'''=\varphi^{(5)}-\varphi^{(4)}>0$. Therefore $\Psi(x)\ge\Psi(0)+\Psi'(0)x+\frac12\Psi''(0)x^2$ for $x\ge0$, and integrating against $d\nu$,
\[
\int\Psi\,d\nu\ \ge\ \Psi(0)+\Psi'(0)+\frac\eps2\,\Psi''(0)=\frac{3(10-\eps)^2(2-\eps)^2}{1024\,\eps}\ \ge\ \frac{3(2-\eps)^2}{16\eps}.
\qedhere
\]
\end{proof}

\section{\texorpdfstring{The chart of $H_i$ (Lemma~\ref{lem:chart})}{The chart of Hi}}\label{app:chart}
\begin{proof}[Proof of Lemma~\ref{lem:chart}]
\emph{The chart.} For $q\in H_i$ the other coordinates are positive with sum $N-t$, so $x\in\mathcal T^\circ$; and $\sum_{j\ne i}q_j^2=N\eps-t^2$ gives $|x|^2=k(t)$. Since
\[
k'(t)=-\frac{2N(t-\eps)}{(N-t)^3}<0\qquad(\eps<t<N),
\]
the condition $t>2\eps$ is equivalent to $|x|^2<k(2\eps)$. Conversely, let $x\in\Omega$. By Cauchy--Schwarz $|x|^2\ge1/(N-1)=k(t^*)$; the last equality holds because at the pole the other coordinates all equal $(N-t^*)/(N-1)$. So there is a unique $t(x)\in(2\eps,t^*]$ with $k(t(x))=|x|^2$, smooth in $x$ because $k'\ne0$, and $F(x)\in H_i$: its coordinates are positive, sum to $t+(N-t)=N$, and have sum of squares $t^2+(N-t)^2k(t)=N\eps$. The two maps are inverse to each other.

\emph{(a) The density.} In the meridian coordinates of Lemma~\ref{lem:slices}(2), $(q_j)_{j\ne i}=c\one+\rho\,\omega$ with $\omega$ in the unit sphere of $\{v\in\R^{N-1}:\sum_jv_j=0\}$, and $d\sigma=\rho^{N-3}\gamma^{-1/2}\,dt\,d\omega$. Dividing by $N-t=(N-1)c$ gives $x=\bar x+r\omega$, with $\bar x=\one/(N-1)$ and $r=\rho/(N-t)$. These are polar coordinates around $\bar x$ in the hyperplane, so $dx=r^{N-3}\,dr\,d\omega$. Since $r^2=|x|^2-\frac1{N-1}=k(t)-\frac1{N-1}$, we have $2r\,dr=k'(t)\,dt$, i.e.\ $|dt/dr|=r(N-t)^3/(N(t-\eps))$. Hence
\[
d\sigma=\frac{\rho^{N-3}}{r^{N-3}}\,\gamma^{-1/2}\Big|\frac{dt}{dr}\Big|\,dx=\frac{(N-t)^{N-1}}{N(t-\eps)}\,\rho\,\gamma^{-1/2}\,dx .
\]
From the definitions of $\rho$ and $\gamma$ (Lemmas~\ref{lem:calc} and~\ref{lem:slices}), $\rho^2=\frac N{N-1}R^2\gamma$, so $\rho\gamma^{-1/2}=R\sqrt{N/(N-1)}=N\sqrt{(\eps-1)/(N-1)}$; this gives $w$. Both factors $(N-t)^{N-1}$ and $1/(t-\eps)$ decrease in $t$. The two-sided bound follows from $-u-u^2\le\log(1-u)\le-u$ for $0\le u\le\frac12$, applied with $u=t/N$ and multiplied by $N-1$.

\emph{(b) The diameter.} $\Omega$ lies in the ball of centre $\bar x$ and radius $r_0$, where $r_0^2=k(2\eps)-\frac1{N-1}<\frac{N\eps}{(N-2\eps)^2}\le\frac{2\eps}N$, because $(N-2\eps)^2\ge N^2/2$ for $N\ge7\eps$. So $\diam(\Omega)^2\le4r_0^2\le8\eps/N$.

\emph{(c) The Lipschitz bound.} Let $v$ be tangent to $\Omega$, i.e.\ $\sum_jv_j=0$. Differentiating $k(t(x))=|x|^2$ gives $k'(t)\,dt(v)=2\langle x,v\rangle=2\langle x-\bar x,v\rangle$, so, using $t-\eps>\eps$, $N-t<N$ and $r_0\le\sqrt{2\eps/N}$,
\[
|dt(v)|\ \le\ \frac{2r_0|v|}{|k'(t)|}=\frac{r_0(N-t)^3}{N(t-\eps)}\,|v|\ \le\ \frac{r_0N^2}{\eps}|v|\ \le\ \sqrt{\tfrac2\eps}\,N^{3/2}|v| .
\]
Since $dF(v)=\big(dt(v),\ ((N-t)v_j-x_j\,dt(v))_{j\ne i}\big)$ and $|x|\le1$,
\[
|dF(v)|^2\ \le\ dt(v)^2+2(N-t)^2|v|^2+2|x|^2dt(v)^2\ \le\ 2N^2|v|^2+3\,dt(v)^2\ \le\ \Big(2+\frac6\eps\Big)N^3|v|^2 .
\]
\end{proof}

\section{Numerical methods}\label{app:numerics}

This appendix describes the numerical methods behind Figures~\ref{fig:sections}, \ref{fig:barrier} and~\ref{fig:lmf_lMC}. None of the proofs uses them.

\subsection{Evaluation of the partition functions}\label{app:numerics-omega}
Figure~\ref{fig:sections} evaluates $\Omega_n(S,E)$ by exponential tilting and an exact decomposition into fluid and condensate. For a cutoff $\ln N\ll L\ll a\sqrt N$, classifying configurations by the number $k$ of sites with $q_i>L$,
\begin{equation}
\Omega_n(S,E)=\sum_{k\ge0}\binom nk\int_{(L,\infty)^k}\!d^kq\;
\Omega^{<L}_{n-k}\Bigl(S-\sum_j q_j,\;E-\sum_j q_j^2\Bigr),
\label{eq:bigjump}
\end{equation}
with $\Omega^{<L}_m$ as in Lemma~\ref{lem:density}. On a bounded support $\cZ_L(\eta,\beta)=\int_0^L e^{-\eta q-\beta q^2}dq$ is finite for every real $\beta$, including $\beta<0$, so the saddle point always exists, and the local limit theorem gives $\ln\Omega^{<L}_m(S',E')=\eta S'+\beta E'+m\ln \cZ_L-\ln\bigl(2\pi m\sqrt{\det C}\bigr)+O(1/m)$, with $C$ the tilted covariance of $(q,q^2)$. The terms $k=1,2$ are computed by quadrature and $k\ge3$ is negligible. Independence of $L$, the normalisation of the one- and two-site marginals, and agreement with an exact-constraint Montecarlo at $N=100$ check the accuracy.

\subsection{Montecarlo dynamics}\label{app:numerics-mc}
We simulated Brownian motion on the microcanonical surface
$M_{N,\eps}=\{q\in\mathbb{R}_{\ge0}^N:\ \sum_i q_i=N,\ \sum_i q_i^2=N\eps\}$ and measured its
relaxation rates directly. We use the technique in \cite{basile2006momentum,iubini2014coarsening}. 

The elementary move picks three distinct sites $(i,j,k)$ at random and
rotates $(q_i,q_j,q_k)$ by an angle $\phi\sim U[-d,d]$ along the circle on which both
$q_i+q_j+q_k$ and $q_i^2+q_j^2+q_k^2$ are fixed. A move that would make a coordinate negative is
rejected. The microcanonical measure restricted to each such circle is uniform in the angle, so
the move satisfies detailed balance exactly for any $d$. Time is calibrated through the exact
identity $\sum_t X_t^2=N R^2\,\Delta_\M$, with $X_t$ the rotation generator of triple $t$ and
$R^2=N(\eps-1)$. Assigning each move the time step
$\delta t=d^2R^2/[3(N-1)(N-2)]$ makes the generator of the simulated process $\Delta_\M$, so the
measured rates are eigenvalues of $-\Delta_\M$, in the same normalisation as $\lambda_{\rm MF}$.
The dynamics commutes with permutations of the sites, so $L^2(\M)$ splits into $S_N$ sectors, each
with its own gap. 

We probed these sectors separately, using time autocorrelations of
\emph{(i)} single-site functions $f(q_i)$, $f\in\{q,q^2,q^3,e^{-q}\}$, from which the combination
$C_{\rm same}-C_{\rm cross}$ isolates the standard sector;
\emph{(ii)} two-site observables
$(\varphi_a(q_i)-\varphi_a(q_l))(\varphi_b(q_j)-\varphi_b(q_k))$ on disjoint quadruples of sites,
which have no trivial or standard component and therefore probe the sectors $(N-2,2)$,
$(N-2,1,1)$, \dots; and \emph{(iii)} collective observables ($\sum_i q_i^3$, $\sum_i e^{-q_i}$,
$\sum_i\sqrt{q_i}$, $\max_i q_i$), which probe the trivial sector.

In each sector the slowest rate was obtained from a generalised eigenvalue problem on the
correlation matrix of the observables, followed by a single-exponential fit on a window selected
from the data alone, with no reference to $\lambda_{\rm MF}$. Pointwise errors on the correlations come from
Bartlett's formula, and errors on the rates from a leave-one-seed-out jackknife over 8
independent seeds. We used $\eps\in\{1.3,1.5,1.7,1.9\}$ and $N\in\{100,\dots,1600\}$, with a
simulated time per seed such that $T\lambda_{\rm MF}\approx3\times10^4$.
Rejection at the faces $q_i=0$ induces a bias in the rates that is linear in the step size, so
every point was run at $d=0.4$, $0.2$ and $0.1$. 

The rates were then extrapolated first linearly
to $d\to0$ and then as $\lambda_{\rm MC}(N)=\lambda_{\rm MC}+c/N$ to $N\to\infty$.
At $\eps=1.9$ we added runs at $N=3200$ and $6400$ with $d=0.05$ and $0.1$. These use reflection
of the angle at the ends of the allowed arc instead of rejection. Reflection is also exactly
reversible and reduces the slope of the $O(d)$ bias by about one third.

\bibliographystyle{apsrev4-2}
\bibliography{references}

\end{document}